\documentclass[a4paper]{amsart}
\usepackage{amsmath,amsthm, amscd, amssymb, amsfonts}
\usepackage[dvips]{graphicx}
\usepackage[all]{xy}
\numberwithin{equation}{section}
\usepackage[english]{babel}
\usepackage{color}

\theoremstyle{plain}

\newtheorem{theorem}{Theorem}[section]
\newtheorem{corollary}[theorem]{Corollary}
\newtheorem{proposition}[theorem]{Proposition}
\newtheorem{lemma}[theorem]{Lemma}

\theoremstyle{definition}

\theoremstyle{remark}
\newtheorem{remark}[theorem]{Remark}

\renewcommand{\1}{\textbf{1}}

\newcommand{\KER}{\mathfrak{Ker}}

\newcommand{\R}{{\mathcal R}}

\newcommand{\C}{{\mathcal C}}

\newcommand\id{\operatorname{id}}
\newcommand\Aut{\operatorname{Aut}}

\newcommand\Tr{\operatorname{Tr}}
\newcommand\gr{\operatorname{gr}}

\newcommand\GL{\operatorname{GL}}

\newcommand\St{\operatorname{St}}

\newcommand\Hom{\operatorname{Hom}}
\newcommand\Rep{\operatorname{Rep}}

\newcommand\Ind{\operatorname{Ind}}

\newcommand\FPdim{\operatorname{FPdim}}

\newcommand\vect{\operatorname{Vec}}

\newcommand\Irr{\operatorname{Irr}}
\newcommand\Pic{\operatorname{G}}
\newcommand\rep{\operatorname{rep}}

\newcommand{\ot}{\otimes}
\newcommand{\mtc}{\mathcal}

\newcommand{\lb}{\label}

\newcommand{\al}{\alpha}

\newcommand{\eps}{\epsilon}
\newcommand{\bn}{\begin}

\newcommand{\cS}{\mathcal{S}}
\newcommand{\cO}{\mtc{O}}

\numberwithin{equation}{section}

\newtheorem{defn}[theorem]{Definition}
\newtheorem{cor}[theorem]{Corollary}
\newtheorem{rem}[theorem]{Remark}
\newtheorem{example}[theorem]{Example}

\newcommand{\bp }{\bn{proposition}}
\newcommand{\ep}{\end{proposition}}

\newcommand{\uw}{\uparrow}

\newcommand{\ch}{\chi}
\newcommand{\mtr}{\mathrm}

\newcommand{\ncm}{\newcommand}\newcommand{\gm}{\gamma}
\numberwithin{equation}{section}
\newcommand{\el}{\end{lemma}}\newcommand{\bl}{\bn{lemma}}
\newcommand{\et}{\end{thm}}
\newcommand{\beqarn}{\begin{eqnarray*}}
\newcommand{\eeqarn}{\end{eqnarray*}}

\newcommand{\eeqn}{\end{equation*}}
\newcommand{\bpf}{\bn{proof}}
\newcommand{\epf}{\end{proof}}
\ncm{\cX}{\mtc{X}}
\ncm{\wt}{\widetilde}
\newcommand{\ra}{\rightarrow}
\ncm{\ro}{\rho}
\ncm{\sg}{\sigma}
\ncm{\np}{\newpage}
\ncm{\ebl}{\end{thebibliography}}
\ncm{\bbl}{\begin{thebibliography}}
\ncm{\chd}{_{ _{\ch}}}
\ncm{\ald}{_{ _{\al}}}
\ncm{\cP}{\mathcal{P}}
\ncm{\ei}{e_i}
\ncm{\eij}{e_{i,\;j}}
\ncm{\bdef}{\begin{defn}}
\ncm{\edf}{\end{defn}}
\ncm{\bc}{\begin{cor}}
\ncm{\ec}{\end{cor}}
\ncm{\er}{\end{rem}}
\ncm{\br}{\begin{rem}}
\ncm{\bd}{\begin{document}}
\ncm{\ed}{\end{document}}
\ncm{\beq}{\begin{equation}}
\ncm{\eeq}{\end{equation}}
\ncm{\bea}{\begin{eqnarray}}
\ncm{\eea}{\end{eqnarray}}
\ncm{\beanon}{\begin{eqnarray*}}
\ncm{\eeanon}{\end{eqnarray*}}\ncm{\ek}{\eps|_K}\ncm{\diez}{\#}
\ncm{\bwt}{\bowtie}
\ncm{\cC}{\mtc{C}}
\ncm{\cs}{\cS}
\ncm{\cA}{\mathcal{A}}
\ncm{\HKer}{\mtr{HKer}}
\ncm{\LKER}{\mtr{LKER}}
\ncm{\aad}{\mtr{ad}}
\ncm{\Dr}{\mtr{D}}
\ncm{\cD}{\mathcal{D}}
\ncm{\Dc}{\mtc{D}}
\ncm{\Vc}{\mtr{Vec}}
\ncm{\cK}{\mtc{K}}
\ncm{\cM}{\mtc{M}}
\ncm{\cE}{\mtc{E}}
\ncm{\chr}{character }
\ncm{\bw}{\bwt}
\ncm{\cY}{\mtc{Y}}
\ncm{\cd}{\cD}
\ncm{\bx}{\boxtimes}
\ncm{\cc}{\cC}\ncm{\irr}{\Irr}
\ncm{\bgm}{\Gamma}
\ncm{\goesto}{\mapsto}
\ncm{\one}{\mathrm{1}_{\cc}}
\ncm{\bvec}{\mtr{Vec}}
\ncm{\omg}{\omega}
\ncm{\og}{\omega}
\ncm{\dt}{\delta}
\ncm{\cg}{\mtc{G}}
\ncm{\blue}{\textcolor[rgb]{.00, .00, 1.00}}
\ncm{\red}{\textcolor[rgb]{1.00, .00, .00}}
\ncm{\green}{\textcolor[rgb]{.00, 1.00, .00}}
\ncm{\dzs}{\#_{\sg}}
\ncm{\md}{\medskip}\ncm{\opl}{\oplus}\ncm{\xra}{\xrightarrow}
\ncm{\Gm}{\Gamma}\ncm{\ind}{\mtr{Ind}}\ncm{\res}{\mtr{Res}}\ncm{\core}{\mtr{core}}\ncm{\lker}{\mtr{LKer}}\ncm{\bb}{\blue}
\renewcommand{\subjclassname}{
\textup{2000} Mathematics Subject Classification}
\ncm{\bq}{\beq}\ncm{\eq}{\eeq}
\ncm{\bne}{\begin{enumerate}}
\ncm{\ene}{\end{enumerate}}
\ncm{\bqarn}{\beqarn}
\ncm{\eqarn}{\eeqarn}
\ncm{\inv}{^{-1}}

\begin{document}

\title[Fusion rules of an equivariantization]{Fusion rules of equivariantizations of fusion categories}
\date{Revised version of October 12, 2019}

\begin{abstract} We determine the fusion rules of the equivariantization of a fusion category $\C$ under the action of a finite group $G$ in terms of the fusion rules of $\C$ and
 group-theoretical data associated to the group action. As an application we obtain a formula for the fusion rules in an equivariantization of a pointed fusion category in terms of group-theoretical data. This entails a description of the fusion rules in any braided group-theoretical fusion category.
\end{abstract}

\author{Sebastian Burciu and Sonia Natale}
\address{S. Burciu: Inst.\ of Math.\ ``Simion Stoilow" of the Romanian Academy
P.O. Box 1-764, RO-014700, Bucharest, Romania and University of Bucharest, Faculty of Mathematics and Computer Science, Algebra and Number Theory Research Center, 14 Academiei St., Bucharest, Romania }\email{sebastian.burciu@imar.ro}
\address{S. Natale: Facultad de Matem\'atica, Astronom\'\i a y F\'\i sica.
Universidad Nacional de C\'ordoba. CIEM-CONICET. (5000) Ciudad
Universitaria. C\'ordoba, Argentina}\email{natale@famaf.unc.edu.ar}
\thanks{The work of S. Burciu was partially supported by a grant of the Romanian National Authority for Scientific Research, CNCS-UEFSCDI, grant no. 88/ 05.10.2011, "Hopf algebras and related topics".The work of S. Natale was partially supported by CONICET and SeCYT--UNC. This project has started while the first author was visiting Universidad Nacional de C\'ordoba. He thanks UNC for its hospitality during his stay.}

\keywords{Fusion category;  fusion rules; Grothendieck ring;  equivariantization}

\subjclass{18D10; 16T05}

\maketitle

\section{Introduction}\label{intro}

Throughout this paper we shall work over an algebraically closed
field $k$ of characteristic zero. Let $\C$ be a fusion category
over $k$, that is, $\C$ is a semisimple rigid monoidal category
over $k$ with finitely many isomorphism classes of simple objects,
finite-dimensional Hom spaces, and such that the unit object $\1$
is simple.

Consider an action $\rho: \underline G \to \underline
\Aut_{\otimes} \C$ of a finite group $G$  by tensor
autoequivalences of $\C$ and let $\C^G$ be the equivariantization
of $\C$ with respect to this action. Equivariantization under a
finite group action, as well as its applications, generalizations
and related constructions, have been intensively studied in the
last years by several authors. See for instance \cite{agaitsgory,
tensor-exact, indp-exact, DGNOI, ENO2, fw, nik, tambara}.

In the sense of the notions introduced in \cite{tensor-exact,
indp-exact}, equivariantization gives rise in a  canonical way to
a \emph{central exact sequence of tensor categories} $$\rep G \to
\C^G \to \C,$$ where $\rep G$ is the category of
finite-dimensional representations of $G$. On the other hand,
combined with the dual notion of (graded) group extension of a
fusion category, equivariantization underlies the notion of
\emph{solvability} of a fusion category developed in \cite{ENO2}.

\medbreak An important invariant of a fusion category $\C$ is its
Grothendieck ring, $\gr(\C)$. For instance, the knowledge of the structure of
the Grothendieck ring allows to determine all fusion subcategories
of $\C$, which correspond to the so-called \emph{based subrings}.

Let $\Irr(\C) = \{ \1 = S_0, \dots, S_n \}$ denote the set of
isomorphism classes of simple objects of $\C$. Then $\Irr(\C)$ is
a basis of $\gr(\C)$ and, for all $0\leq i, j
\leq n$, we have a relation
\begin{equation}\label{nijl} S_iS_j = \sum_{l = 0}^n N_{i, j}^l \,
S_l, \end{equation} where $N_{i, j}^l$ are non-negative integers
given by $N_{i, j}^l = \dim \Hom_\C(S_l,  S_i \otimes S_j)$, $0
\leq l \leq n$. The relations \eqref{nijl} are known as the
\emph{fusion rules} of $\C$. They are determined
by the set $\Irr(\C)$ and the multiplicities $N_{i, j}^l \in
\mathbb Z_{\geq 0}$.

The main result of this paper is the determination of the fusion
rules of $\C^G$ in terms of the fusion rules of $\C$ and certain
canonical group-theoretical data associated to the group action.
This is contained in Theorem \ref{tp}. As it turns out, the
structure of the Grothendieck ring of $\cc^G$ resembles that of
the rings introduced by Witherspoon in \cite{ws}. See Remark \ref{rings}.

As an example, consider a semisimple cocentral Hopf
algebra extension $H$ of a Hopf algebra $A$ by a finite group $G$,
that is, $H$ fits into a cocentral exact sequence
$$k \to A \to H \to kG \to k.$$
As shown in \cite{ext-ty} the category $\Rep H$ of
finite-dimensional representations of $H$ is an equivariantization
$(\Rep A)^G$ with respect to an appropriate action of $G$ on $\Rep
A$. Thus Theorem \ref{tp} implies that the fusion rules of the
category $\Rep H$ can be described in terms of the fusion rules of
$\Rep A$ and the action of $G$. In particular, Theorem \ref{tp}
generalizes the results obtained for cocentral abelian extensions
of Hopf algebras in \cite{goff} and \cite{ws}.

\medbreak We discuss in detail the case where $\C$ is a
\emph{pointed} fusion category, that is, when all simple objects
of $\C$ are invertible. In this case the fusion rules of $\C^G$
are described completely in terms of group-theoretical data. See
Theorem \ref{tp-pointed}.

It is known that every braided group-theoretical fusion category
is an equivariantization of a pointed fusion category \cite{NNW,
naidu}. Therefore, our results entail the determination of the
fusion rules in any braided group-theoretical fusion category.

\medbreak In order to establish  Theorem \ref{tp}, we give an
explicit description of the simple objects of the
equivariantization $\cc^G$. 
Such an action induces naturally an action of $G$ on the set
$\Irr (\C)$ of isomorphism classes of simple objects of $\C$. Let
$Y \in \Irr(\C)$ and let $G_{Y} \subseteq G$ denote the inertia
subgroup of $Y$, that is, \bq G_{Y}=\{g \in G\;|\; \ro^g(Y)\simeq
Y\}. \eq

We show that isomorphism  classes of simple objects of $\cc^G$ are
parame\-terized by pairs $(Y, \pi)$, where $Y$ runs over a set  $\Irr(\C)/G = \{Y_1, \dots, Y_s\}$ of simple objects of $\C$ whose isomorphism classes represent the orbits of $G$ in $\Irr(\C)$, and $\pi$ runs over the isomorphism classes of irreducible projective representation of
the inertia subgroup $G_{Y}$ with a certain factor set $\widetilde \alpha_Y
\in Z^2(G_Y, k^*)$. This result is analogous to the
parameterization of  irreducible representations of a finite group
in terms of those of a normal subgroup given by Clifford Theorem.
It extends the description obtained in \cite{MW} for the case
where $\C^G$ is the category of representations of an (algebra)
group crossed product (see \cite[Subsection 3.1]{ext-ty}).

In the case where $\C$ is a fusion category, the duality in $\C$ gives rise to a ring involution ${}^*:\gr(\C) \to \gr(\C)$.
We describe this ring involution for the Grothendieck ring of $\C^G$ in Subsection \ref{involution}, more precisely, we use Theorem \ref{tp} in order to determine  the dual of the simple object in $\C^G$ corresponding to a pair $(Y, \pi)$ as above.

We may regard $\rep G$ as a fusion subcategory of $\C^G$ under a
canonical embedding. In this way $\C^G$ becomes a $\rep
G$-bimodule category under the action given by tensor product. As
another consequence of Theorem \ref{tp}, we give a decomposition
of $\C^G$ into indecomposable $\rep G$-module categories. See
Theorem \ref{bimod-equiv}.

\medbreak The paper is organized as follows. In Section
\ref{simples-equiv} we recall the definition of the
equivariantization of a semisimple abelian category over $k$ under
a finite group action and give a parameterization of its simple
objects. With respect to this parameterization, we determine in
Section \ref{fusion-rules} the fusion rules in an
equivariantization of a fusion category. Appart from the main
result, Theorem \ref{tp},  we also present in this section the
above mentioned applications to the determination of the dual of a
simple object and the decomposition of $\C^G$ as a $\rep
G$-bimodule category. In Section \ref{pointed-fusion} we
specialize our main result to the case of an equivariantization of
a pointed fusion category and in particular, to braided
group-theoretical fusion categories. We include an Appendix at the
end of the paper, where we give an account of the relevant facts
about projective group representations needed throughout.

%

\section{Simple objects of an equivariantization}\label{simples-equiv}

The goal of this section is to describe a (mostly well-known)
Clifford correspondence entailing a classification of isomorphism
classes of simple objects of $\C^G$ in terms of the action of $G$
on the set $\Irr(\C)$ of isomorphism classes of simple objects of
$\C$. 
Some instances of this correspondence appear for instance in \cite[Proof of Proposition 6.2]{ENO2} and \cite[Proposition 5.5]{naidu}.

By abuse of notation, we often indicate an object of a
category $\C$ and its isomorphism class by the same letter.

\subsection{Equivariantization under a finite group action}

Let $\C$ be a finite semisimple $k$-linear
category and let $G$ be a finite group. Let also $\rho: \underline
G \to \underline \Aut \, \C$ be an action of $G$ on $\C$ by
$k$-linear autoequivalences. Thus, for every $g \in G$, we have a
$k$-linear functor $\rho^g: \C \to \C$ and natural isomorphisms
$$\rho^{g,h}_2 : \rho^g \rho^h \to \rho^{gh}, \quad g, h \in G,$$ and
$\rho_0 : \id_\C \to \rho^e$, subject to the following conditions:
\begin{align}\label{ro-2} & (\rho^{ab, c}_2)_X \, (\rho^{a,b}_2)_{\rho^c(X)} =
(\rho^{a, bc}_2)_X \, \rho^a((\rho^{b,c}_2)_X), \\ & (\rho^{a,
e}_2)_X \rho^{a}_{\rho_0(X)} = (\rho^{e, a}_2)_X
(\rho_0)_{\rho^{a}(X)},
\end{align}
for all objects $X \in \C$, and for all $a, b, c \in G$.
By the naturality of $\rho^{g,h}_2$, $g, h \in G$, we have the following relation:
\begin{equation}\label{nat-ro}\rho^{gh}(f) \, (\rho_2^{g, h})_{Y} = (\rho_2^{g, h})_X \, \rho^g\rho^h(f),
\end{equation} for every morphism $f: Y \to X$ in $\C$.
For simplicity, we shall assume in what follows that $\rho^e = \id_\C$ and
$\rho_0$, $\rho^{g, e}_2$, $\rho^{e, g}_2$ are identities.

\medbreak Let $\C^G$ denote the corresponding
\emph{equivariantization}. Recall that $\C^G$ is a finite
semisimple $k$-linear category whose objects are $G$-equivariant
objects of $\C$, that is, pairs $(X, \mu)$, where $X$ is an object
of $\C$ and $\mu = (\mu^g)_{g \in G}$, such that $\mu^g:\rho^gX
\to X$ is an isomorphism, for all $g \in G$, satisfying
\begin{equation}\label{deltau} \mu^g \rho^g(\mu^h) = \mu^{gh} (\rho^{g,
h}_2)_X, \quad \forall g, h \in G, \quad \quad
\mu_e{\rho_0}_X=\id_X.\end{equation} A morphism $f: (X, \mu) \to
(X', \mu')$ in $\C^G$ is a morphism $f: X \to X'$ in $\C$ such
that $f\mu^g = {\mu'}^g\rho^g(f)$, for all $g \in G$.

\medbreak We shall also say that an object $X$ of $\C$ is
\emph{$G$-equivariant} if there exists such a collection $\mu =
(\mu^g)_{g \in G}$ so that $(X, \mu) \in \C^G$. Note that $\mu$ is
not necessarily unique.

\medbreak The forgetful functor $F: \C^G \to \C$, $F(X, \mu) = X$,
is a dominant functor. The functor $F$ has a left adjoint $L: \C
\to \C^G$, defined by $L(X) = (\oplus_{t \in G}\rho^tX, \mu_X)$,
where $(\mu_X)_g: \oplus_{t \in G}\rho^g\rho^tX \to \oplus_{t \in
G}\rho^tX$ is given componentwise by the isomorphisms
$(\rho^{g,t}_2)_X$. The composition $T_\rho = FL: \C \to \C$ is a
faithful $k$-linear monad on $\C$ such that $\C^G$ is equivalent to the
category $\C^{T_\rho}$ of $T_\rho$-modules in $\C$. See
\cite[Subsection 5.3]{tensor-exact}.

\medbreak 
We assume in what follows that $\C$ is a fusion category
over $k$ and $\rho: \underline G \to \underline \Aut_{\otimes} \C$
is an action of $G$ on $\C$ by \emph{tensor autoequivalences},
that is, $\rho_2^{g, h}: \rho^g \rho^h \to \rho^{gh}$ are natural
isomorphisms of tensor functors, for all $g, h \in G$. Thus, for
all $g \in G$, $\rho^g$ is endowed with a monoidal structure
$(\rho_2^g)_{X, Y}: \rho^g(X \otimes Y) \to \rho^g(X) \otimes
\rho^g(Y)$, $X, Y \in \C$, and the following relation holds:
\begin{equation}\label{tensor-rho} {\rho_2^{gh}}_{X, Y} {\rho_2^{g, h}}_{X \otimes Y} = ({\rho_2^{g, h}}_{X} \otimes {\rho_2^{g, h}}_{Y}) \, {\rho_2^{g}}_{\rho^hX, \rho^hY} \, {\rho^g(\rho_2^{h}}_{X, Y}),
\end{equation}
for all $g, h \in G$, $X, Y \in \C$.

\medbreak Then $\C^G$ is also a fusion category with tensor product
$(X, \mu_X) \otimes (Y, \mu_Y) = (X \otimes Y, (\mu_X \otimes \mu_Y){\rho_2}_{X, Y})$, where for all $g \in G$, ${\rho_2}_{X, Y}^g: \rho^g(X \otimes Y) \to \rho^g(X) \otimes \rho^g(Y)$ is
the monoidal structure on $\rho^g$. In particular the forgetful functor $F: \C^G \to \C$ is a strict monoidal functor.

\subsection{Frobenius-Perron dimensions of simple objects of $\C^G$} Let $X, Y \in \C$. Then $\Hom_\C(\rho^gX, \rho^gY)
\simeq \Hom_\C(X, Y)$, for all $g \in G$. Therefore, for all $g
\in G$, and for all objects $M$ of $\C^G$, we have
\begin{equation}\label{multY}\Hom_\C(F(M), \rho^gY) \simeq
\Hom_\C(F(M), Y). \end{equation}

The action of the group $G$ on $\C$ permutes isomorphism classes
of simple objects of $\C$. Let $Y \in \Irr(\C)$. We shall denote
$G_Y: = \St_G(Y) \subseteq G$ the \emph{inertia} subgroup of $Y$, that is,
\begin{equation*} G_{Y} = \{ g \in G|\; \rho^g(Y) \simeq Y\},
\end{equation*}  Then $Y$ has exactly $n = [G: G_Y]$ mutually
nonisomorphic $G$-conjugates $Y = Y_1, \dots, Y_n$. For every $1
\leq j \leq n$, we have $Y_j \simeq \rho^{g_j}Y$, where $g_1 = e,
\dots, g_n$ is a complete set of representatives of the left
cosets of $G_Y$ in $G$.

\begin{proposition}\label{conjugates} Let $M = (X, \mu)$ be a simple object of $\C^G$ and let
$Y$ be a simple constituent of $X$ in $\C$. Let $Y = Y_1, \dots,
Y_n$, $n = [G: G_Y]$, be the mutually nonisomorphic $G$-conjugates
of $Y$. Then $X \simeq m \oplus_{i = 1}^nY_i$, where $m = \dim
\Hom_\C(X, Y)$.
\end{proposition}

\begin{proof} Consider the object $T(Y) = FL(Y) = \oplus_{g \in G}\rho^gY$. Let $Z$ be a simple object of $\C$ such that $Z \ncong Y_j$, $j = 1, \dots, n$. Then $\Hom_\C(Z, T(Y)) = 0$.
By adjointness, we have $\Hom_{\C^G}(L(Y), M) \simeq \Hom_\C(Y, X)
\neq 0$. Then $M$ is a simple direct summand of $L(Y)$ in $\C^G$.
This implies that $X = F(M)$ is a direct sum of simple subobjects
of $FL(Y) = T(Y)$. Therefore $\Hom_\C(Z, X) = 0$.

Hence $X \simeq \oplus_{i = 1}^nm_i Y_i$, where $m_i = \dim
\Hom_\C(Y_i, X)$, for all $i$. By \eqref{multY}, we have $m_i =
m_1 = m$, for all $i = 1, \dots, n$. This proves the proposition.
\end{proof}

\begin{corollary} Let $M = (X, \mu)$ be a simple object of $\C^G$ and let $Y$ be a simple constituent of $X$ in $\C$. Then $\FPdim M = m [G:G_Y] \FPdim Y$, where $m = \dim \Hom_\C(Y, X)$. \qed \end{corollary}

\subsection{Equivariantization and projective group representations} Let $Y \in \C$ be a fixed simple object. The action $\rho$ of $G$
on $\C$ induces by restriction  an action of $G_Y$ on $\C$ by
autoequivalences. We may thus consider the equivariantization
$\C^{G_Y}$.

By definition of $G_Y$, there exist isomorphisms $c^g:\rho^g(Y)
\to Y$, for all $g \in G_Y$. For all $g, h \in G_Y$, the
composition $c^g \rho^g(c^h)({\rho^{g,
h}_{2_Y}})^{-1}(c^{gh})^{-1}$  defines an isomorphism $Y \to Y$.
Since $Y$ is a simple object, there exist  nonzero
$\widetilde \alpha_Y(g, h) \in k$ such that
\begin{equation}\label{alfa}  \widetilde\alpha_Y(g, h)^{-1} \id_Y = c^g \rho^g(c^h)({\rho^{g, h}_{2_Y}})^{-1}(c^{gh})^{-1}: Y \to Y. \end{equation}
This defines a map $\widetilde\alpha_Y: G_Y \times G_Y \to k^*$
which is a $2$-cocycle on $G_Y$.

\begin{remark}\label{dependence-alfa} The cocycle $\widetilde\alpha_Y$ measures the possible obstruction for $(Y, c)$ to be a $G_Y$-equivariant object, where $c = (c^g)_{g \in G_Y}$.

Consider another choice of isomorphisms $v^g: \rho^g(Y) \to Y$, $g
\in G_Y$. Since $Y$ is a simple object, the composition
$c^g(v^g)^{-1}: Y \to Y$ is given by scalar multiplication by some
$f(g) \in k^*$, for all $g \in G_Y$.  Denoting by
$\widetilde\beta_Y$ the 2-cocycle related to $(v^g)_g$, it easy to
see that $\widetilde\alpha_Y$ and $\widetilde\beta_Y$ differ by
the coboundary of the cochain $f: G_Y \to k^*$. This implies that
the cohomology class $\alpha_Y \in H^2(G_Y, k^*)$ of
$\widetilde\alpha_Y$ depends only on the isomorphism class of the
simple object $Y$.
\end{remark}

\begin{lemma}\label{proj-rep} Let $(X, \mu) \in \C^G$ and let $Y \in \Irr(\C)$.
Consider, for every $g \in G_Y$, isomorphisms $c^g:\rho^g(Y) \to
Y$ and let $\widetilde\alpha_Y$ be the associated 2-cocycle on
$G_Y$. Then the space $\Hom_\C(Y, X)$ carries a projective
representation of $G_Y$ with factor set $\widetilde\alpha_Y$,
defined in the form
\begin{equation}\label{pr}\pi (g) (f) = \mu^g \rho^g(f) (c^g)^{-1} : Y \to X,\end{equation}
for all $f \in \Hom_\C(Y, X)$.
\end{lemma}

\begin{proof} The relation $\pi(g)\pi(h) = \widetilde\alpha_Y(g, h) \pi(gh)$, $g, h \in G$, follows by straightforward computation, using the compatibility conditions for $\rho$. \end{proof}

\begin{remark} Suppose that $\phi: (X, \mu) \to (X', \mu')$
is an isomorphism in $\C^G$. Then the induced isomorphism
$\Hom_\C(Y, X) \to \Hom_\C(Y, X')$, $f \mapsto \phi f$, is an
isomorphism of projective representations. Similarly, if $Y'
\simeq Y$ is another representative of the  isomorphism class of
$Y$ and ${c'}^g:\rho^g(Y') \to Y'$, $g \in G_Y$, is a collection
of isomorphisms, then the projective representations $\Hom_\C(Y,
X)$ and $\Hom_\C(Y', X)$ are projectively equivalent.
\end{remark}

\begin{proposition}\label{cociclo} Let $Y$ be a simple object of $\C$. Let also $c^g:\rho^g(Y) \to Y$,  $g \in G_Y$, be a fixed choice of isomorphisms, and let $\widetilde\alpha_Y$ be the associated 2-cocycle. There is a bijective correspondence between isomorphism classes of simple objects $L = (N, \nu)$ of $\C^{G_Y}$ such that $N$ contains $Y$ as a simple constituent  and isomorphism classes of irreducible  $\widetilde\alpha_Y$-projective representations of the group $G_Y$. If the simple object $L = (N, \nu)$  corresponds to the projective representation $\pi$, then $L \cong \Hom_\C(Y, N) \otimes Y$ and $\pi$ is isomorphic to  $\Hom_\C(Y, N)$. In particular $\FPdim L = \dim \pi \, \FPdim Y$.\end{proposition}

\begin{proof}  Let also $\pi$ be a projective representation of $G_Y$  on the vector space $V$ with factor set $\widetilde\alpha_Y$. Then the pair $(V \otimes Y, \nu)$ is a $G_Y$-equivariant object, where \bq\label{yip}\nu^g = \pi(g) \otimes c^g: \rho^g(V \otimes Y) = V \otimes \rho^g( Y) \to V \otimes Y.\eq 

Conversely, if $L = (N, \nu)$ is an object of $\C^{G_Y}$ such that $N$ contains $Y$ as a simple constituent, then $V = \Hom_\C(Y, N)$ carries a projective representation $\pi$ of $G_Y$ with factor set $\widetilde\alpha_Y$ defined by \eqref{pr}.

These assignments are functorial and mutually inverse up to isomorphisms  on simple objects. This implies the proposition.   \end{proof}

\subsection{The relative adjoint} Consider the forgetful functor $F_Y: \C^G \to \C^{G_Y}$.
We discuss in this subsection the left adjoint $L_Y: \C^{G_Y} \to \C^G$ of the functor $F_Y$.

\medbreak Let $\R$ be a set of representatives of the left cosets
of $G_Y$ in $G$. So that $G$ is a disjoint union $G = \cup_{t \in
\R}tG_Y$.

\medbreak Set, for all $(N, \nu) \in \C^{G_Y}$, $L_Y(N, \nu) =
L_Y^\R(N, \nu) = (\oplus_{t \in \R}\rho^t(N), \mu)$, where, for
all $g \in G$, $\mu^{g}: \oplus_{t \in \R}\rho^g\rho^t(N) \to
\oplus_{t \in \R}\rho^t(N)$ is defined componentwise by the
formula \begin{equation}\label{mu}\mu^{g, t} = \rho^s(\nu^h)
(\rho_2^{s, h})^{-1} \rho_2^{g, t}: \rho^g \rho^t(N) \to
\rho^s(N),\end{equation} where the elements $h \in G_Y$, $s \in \R$
are uniquely determined by the relation
\bq\label{indexch}
gt = sh.
\eq

\begin{remark}\label{indep} We shall show in Proposition \ref{adjointness}
below that the functor $L_Y^\R$ thus defined is left adjoint of
the functor $F_Y$.
By uniqueness of the adjoint, it will follow that, up to isomorphism,
$L_Y^\R$ does not depend on the particular choice of the set
of  representatives $\R$. \end{remark}

\begin{lemma}\label{buena-def} Let $(N, \nu) \in \C^{G_Y}$. Then $L_Y(N, \nu) \in
\C^G$. \end{lemma}

\begin{proof} For every $g \in G$, $t \in \R$, let $s(g, t) \in
\R$, $h(g, t) \in G_Y$ be the elements uniquely determined by the
relation $gt = s(g, t)h(g, t)$. Note that, for
all $a, b \in G$, $t \in \R$, the following relations hold:
\begin{align}\label{s-h}& s(ab, t) = s(a, s(b, t)), \\ \label{habt}& h(ab, t) = h(a, s(b,
t))h(b, t),\\ \label{s-h2}& s(ab, t) h(a, s(b, t)) = a s(b, t).
\end{align}

In order to prove the lemma we shall show that, for all objects $(N, \nu) \in \C^{G_Y}$, and for
all $a, b \in G$, $t \in \R$, the following diagram is
commutative:
\begin{equation}\label{diagrama}\xymatrix@C5em{ \rho^a \rho^b
\rho^t(N) \ar[r]^{(\rho^{a,b}_2)_{\rho^t(N)}} \ar[d]_{\rho^a
(\mu^{b,t})} & \rho^{ab}\rho^t(N)
\ar[d]^{\mu^{ab, t}}\\
\rho^{a}\rho^{s(b, t)}(N)\ar[r]_{\mu^{a, s(b, t)}} & \rho^{s(ab,
t)}(N).}\end{equation}
This is done as follows.
By \eqref{mu}, the relevant maps in diagram \eqref{diagrama} are
given by
\begin{align} \mu^{ab, t}  & = \rho^{s(ab, t)}(\nu^{h(ab, t)}) \,
(\rho_2^{s(ab, t), h(ab, t)})^{-1} \, \rho_2^{ab, t}, \\
\mu^{a, s(b, t)} & = \rho^{s(a, s(b, t))}(\nu^{h(a, s(b, t))}) \,
(\rho_2^{s(a, s(b, t)), h(a, s(b, t))})^{-1} \rho_2^{a, s(b, t)},
\\
\label{form}\rho^a(\mu^{b, t}) & = \rho^a\rho^{s(b, t)} (\nu^{h(b, t)}) \,
\rho^a(\rho_2^{s(b, t)h(b, t)})^{-1} \, \rho^a(\rho_2^{b, t}).
\end{align}
Using \eqref{s-h} and the fact that $(N, \nu)$ is
$G_Y$-equivariant, we compute
\begin{align*}  \rho^{s(ab, t)}(\nu^{h(ab, t)})
& = \rho^{s(a, s(b, t))}(\nu^{h(a, s(b, t)) h(b, t)}) \\
& = \rho^{s(a, s(b, t))} \left( \nu^{h(a, s(b, t))} \, \rho^{h(a,
s(b, t))}(\nu^{h(b, t)}) \, (\rho_2^{h(a, s(b, t)), h(b, t)})^{-1} \right) \\
& = \rho^{s(a, s(b, t))} \left( \nu^{h(a, s(b, t))} \right)
\rho^{s(a, s(b, t))} \rho^{h(a, s(b, t))}\left( \nu^{h(b, t)}
\right) \\ & \qquad \qquad \qquad \times \rho^{s(a, s(b, t))}
\left( \rho_2^{h(a, s(b, t)), h(b, t)}\right)^{-1}
\end{align*}
Similarly, relations \eqref{habt} and \eqref{s-h2}, together with the defining
condition \eqref{ro-2} on the isomorphisms $\rho_2$, give
\begin{align*}  \left(\rho_2^{s(ab, t), h(ab, t)}\right)^{-1}
& =  \rho^{s(ab, t)}\left( \rho_2^{h(a, s(b, t)), h(b, t)} \right)
\left(\rho_2^{s(ab, t), h(a, s(b,t))}\right)_{\rho^{h(b,
t)}(N)}^{-1}  \\ & \times \left( \rho_2^{as(b, t),
h(b, t)} \right)^{-1}.
\end{align*}
The naturality condition \eqref{nat-ro} on $\rho_2^{g, h}$ and relation \eqref{s-h2} imply
that
\begin{align*} \left(\rho_2^{s(ab, t), h(a, s(b,t))}\right)_{\!\rho^{h(b,
t)}(N)}^{-1} = & \left( \rho^{s(ab, t)}\rho^{h(a, s(b,
t))}(\nu^{h(b, t)}) \right)^{\!-1} \! \left( \rho_2^{s(ab, t),
h(a, s(b,t))}\right)^{\!-1} \\ & \times \rho^{s(ab, t), h(a,
s(b,t))}(\nu^{h(b, t)}).
\end{align*}
Hence we get
\begin{align*}  \rho^{s(ab, t)}(\nu^{h(ab, t)})  \left(\rho_2^{s(ab, t), h(ab,
t)}\right)^{\!-1} & = \rho^{s(a, s(b, t))} \left( \nu^{h(a, s(b,
t))} \right) \left( \rho_2^{s(ab, t), h(a, s(b,t))}\right)^{\!-1}
\\ & \; \times \rho^{as(b, t)}(\nu^{h(b, t)}) \left( \rho_2^{as(b, t)),
h(b, t)} \right)^{-1}.
\end{align*}
Composing this resulting morphism with the inverse of $$\rho^{s(a,
s(b, t))}(\nu^{h(a, s(b, t))}) \, (\rho_2^{s(a, s(b, t)), h(a,
s(b, t))})^{-1},$$ and using \eqref{s-h2}, we obtain the
expression
\begin{equation} \rho^{a s(b, t)}(\nu^{h(b, t)}) \left( \rho_2^{as(b, t), h(b, t)}\right)^{-1}.\end{equation}
Using relation \eqref{ro-2}, we see that commutativity of the diagram
\eqref{diagrama} is equivalent to
\begin{align}\label{equivalent}\rho^{a s(b, t)}(\nu^{h(b, t)}) \left( \rho_2^{as(b, t), h(b, t)}\right)^{-1} \rho_2^{a, bt} \, \rho^a(\rho_2^{b, t})
= \mu^{a, s(b, t)} \, \rho^a(\mu^{b, t}).
\end{align}
Finally, we compute
\begin{align*}\rho_2^{a, bt} & = \rho_2^{a, s(b, t)h(b, t)} \\
& = \rho_2^{as(b, t), h(b, t)} (\rho_2^{a, s(b,t)})_{\rho^{h(b,
t)}(N)} \, \rho^a(\rho_2^{s(b, t), h(b, t)})^{-1} \\
& = \rho_2^{as(b, t), h(b, t)} \rho^{as(b, t)}(\nu^{h(b, t)})^{-1}
\rho_2^{a, s(b, t)} \, \rho^a \rho^{s(b, t)} (\nu^{h(b, t)}) \,
\rho^a (\rho_2^{s(b, t), h(b,t)})^{-1}. \end{align*} Combining
this with \eqref{form} we get relation \eqref{equivalent}. This shows that the
diagram \eqref{diagrama} is commutative, as claimed, and finishes
the proof of the lemma. \end{proof}

In view of Lemma \ref{buena-def} there is a well defined functor
$L_Y = L_Y^\R: \C^{G_Y} \to \C^G$.

\begin{proposition}\label{adjointness} The functor $L_Y^\R$ is left adjoint of the
functor $F_Y$. \end{proposition}

\begin{proof} We define natural transformations $\eta\!: \id_{\C^{G_Y}}\! \to
\!F_YL_Y^\R$ and $\epsilon\!: L_Y^\R F_Y \! \to \id_{\C^G}$, in
the form
\begin{align*}& \eta_{(N, \nu)} = i_e: N = \rho^e(N) \to F_YL_Y^\R(N, \nu) =
F_Y(\oplus_{t \in \R}\rho^t(N), \mu),\\
& \epsilon_{(M, u)} = \oplus_{t \in \R}u^t: L_Y^\R F_Y(M, u) =
(\oplus_{t \in \R}\rho^t(M), \mu) \to (M, u), \end{align*} for
every $(N, \nu) \in \C^{G_Y}$, $(M, u) \in  \C^G$.
It is straightforward to verify that $\eta$ and $\epsilon$ are
well-defined and that the compositions
$$F_Y \overset{\eta F_Y}\to F_YL_Y^\R F_Y \overset{F_Y\epsilon}\to F_Y,
\quad L_Y^\R \overset{L_Y^\R \eta}\to L_Y^\R F_YL_Y^\R
\overset{\epsilon L_Y^\R}\to L_Y^\R$$ are identities. This implies
the proposition. \end{proof}

\begin{remark} Note that the restriction of $L_Y$ to the fusion subcategory $\rep G_Y$ of $\C^{G_Y}$ is isomorphic to the induction functor $\rep G_Y \to \rep G \subseteq \C^G$.
\end{remark}

As pointed out in Remark \ref{indep}, we have the following:

\begin{corollary} The functor $L_Y^\R$ is, up to isomorphism, uniquely determined by the subgroup $G_Y$.  \qed \end{corollary}

\subsection{Parameterization of simple objects} The following theorem is the main result of this
section.

\begin{theorem}\label{correspondence} Let $Y \in \C$ be a simple object. Then the functor $L_Y: \C^{G_Y} \to \C^G$ induces a bijective correspondence between isomorphism classes of:

\textrm{(a)} Simple objects $(N, \nu) \in \C^{G_Y}$ such that
$\Hom_\C(Y, N) \neq 0$, and

\textrm{(b)} Simple objects $(X, \mu) \in \C^G$ such that
$\Hom_\C(Y, X) \neq 0$.

If $(X, \mu)$ in $\C^G$ as in (b) corresponds to $(N, \nu)$ in
$\C^{G_Y}$ as in (a), then we have $\Hom_\C(Y, N) \simeq
\Hom_\C(Y, X)$ as projective representations of $G_Y$. Moreover, $N \simeq \Hom_\C(Y, N) \otimes Y$.
\end{theorem}

\begin{proof} Let $(N, \nu) \in \C^{G_Y}$ be a simple object as in (a).
 By Proposition \ref{conjugates} applied to $G_Y$,  $N \simeq m Y$, where $m = \dim \Hom_\C(Y, N)$. Thus $N \simeq \Hom_\C(Y, N) \otimes Y$.

Let $(X, \mu) \in \C^G$ be a simple object such that $(N, \nu)$ is a
simple direct summand of $F_Y(X, \mu)$ in $\C^{G_Y}$. By
adjunction, $(X, \mu)$ is a simple direct summand of $L_Y(N, \nu)$
in $\C^{G}$. Then $X$ is a direct summand of $\oplus_{t \in
G/G_Y}\rho^t(N)$ in $\C$. Since $X$ is $G$-equivariant, then
$\Hom_\C(X, N) \neq 0$.
\\
Therefore $\Hom_\C(Y, X) \neq 0$ and $(X, \mu)$ satisfies the
condition in (b).

Again by Proposition \ref{conjugates}, we get
that  $X \simeq e \oplus_{i = 1}^nY_i$, where  $Y = Y_1, \dots,
Y_n$, $n = [G: G_Y]$, are the mutually nonisomorphic
$G$-conjugates of $Y$. Note that $e = \dim \Hom_\C(Y, X) \leq \dim
\Hom_\C(Y, N) = m$, because the multiplicity of $Y$ in $\rho^t(N)$
is $0$, for all $t \notin G_Y$.

\medbreak Since $(N, \nu)$ is a direct summand of $F_Y(X, \mu)$ in
$\C^{G_Y}$, comparing Frobe\-nius-Perron dimensions, we obtain
$$[G:G_Y] m \FPdim Y = \FPdim L_Y(N, \nu) \leq \FPdim (X, \mu) = e [G:G_Y] \FPdim Y.$$ Therefore $e = m$, and necessarily $(X, \mu) = L_Y(N, \nu)$.
This implies surjectivity of the map induced by $L_Y$ from  (a) to
(b).

Suppose $(N', \nu') \ncong (N, \nu)$ is a simple summand of
$F_Y(X, \mu)$, with $(N', \nu')$ as in (b). Applying the forgetful functor $\C^{G_Y} \to \C$
and comparing the multiplicity of $Y$ we get
$$e = \dim \Hom_\C(Y, X) \geq \dim \Hom_\C(Y, N \oplus N') = m + \dim \Hom_\C(Y, N') > m.$$
This is a contradiction since $e = m$. Hence $(N, \nu)$ is the
unique, up to isomorphisms, simple object as in (b) such that
$$\Hom_{\C^G}(L_Y(N, \nu), (X, \mu)) \simeq \Hom_{\C^{G_Y}}((N,
\nu), F_Y(X, \mu)) \neq 0.$$ This proves injectivity of the map
induced by $L_Y$. Thus this map is bijective, as claimed.

Now suppose that the class of the simple object $(X, \mu)$ of $\C^G$ as in (b) corresponds to the class of the simple object $(N, \nu)$ of
$\C^{G_Y}$ as in (a). The proof above shows that $N \simeq \Hom_\C(Y, N) \otimes Y$.

\medbreak As we have shown,  $\dim \Hom_\C(Y, X) = \dim \Hom_\C(Y, N)$. Since $N$ is a direct summand of $X$ in $\C$, then
$\Hom_\C(Y, N) \subseteq \Hom_\C(Y, X)$, thus these spaces are equal. Furthermore, the corresponding projective representations given by Formula \eqref{pr} clearly coincide on both spaces. This finishes the proof of the theorem.
\end{proof}

Let $\{Y_1, \dots, Y_s\}$ be a set of simple objects of $\C$ whose isomorphism classes represent the orbits of $G$ in $\Irr(\C)$. By abuse of notation, the set $\{Y_1, \dots, Y_s\}$ will be indicated by $\Irr(\C)/G$ in what follows. 

\medbreak
Combining Theorem \ref{correspondence} with Proposition
\ref{cociclo} we obtain the following:

\begin{corollary}\label{cor-clifford} There is a bijective correspondence between the set of isomorphism classes of simple objects of $\C^G$
and the set of pairs $(Y, \pi)$, where $Y$ runs over the  set $\Irr(\C)/G = \{Y_1, \dots, Y_s\}$ and $\pi$ runs over the
isomorphism classes of irreducible $\widetilde \alpha_Y$-projective
representations  of the inertia subgroup $G_Y \subseteq G$, where
 $\widetilde \alpha_Y$ is the cocycle determined by \eqref{alfa}.

Let $(X, \mu)$ be the simple object corresponding to the pair $(Y,
\pi)$. Then we have $X \simeq \oplus_{t \in G/G_Y}\rho^t(V_\pi
\otimes Y)$. In particular, $\FPdim (X, \mu) = \dim \pi [G:G_Y]
\FPdim Y$. \qed \end{corollary}

\begin{proof} Let $Y \in \Irr(\C)/G$ and let $\pi$ be an irreducible $\widetilde \alpha_Y$-projective representation of $G_Y \subseteq G$. By Proposition \ref{cociclo}, $V_\pi \otimes Y$ is a simple object of $\C^{G_{Y}}$ and moreover $V_\pi \otimes Y \cong V_{\pi'} \otimes Y$ if and only if $\pi \cong \pi'$. By Theorem \ref{correspondence},  $L_Y(V_\pi \otimes Y)$ is a simple object of $\C^G$ and $L_Y(V_\pi \otimes Y) \cong L_Y(V_{\pi'} \otimes Y)$ if and only if $\pi \cong \pi'$.
In addition, the simple objects $L_Y(V_\pi \otimes Y)$ and $L_{Y'}(V_{\pi'} \otimes Y')$ are not isomorphic if $Y \neq Y' \in \Irr(\C)/G$. 

Finally note that $\FPdim L_Y(V_\pi \otimes Y) = \dim \pi [G:G_Y]  \FPdim Y$. Hence, as $Y$ runs over $\Irr(\C)/G$ and $\pi$ runs over a set  $\Irr_{\widetilde{\alpha}_{Y}}$ of representatives of the isomorphism classes of irreducible projective $\widetilde{\alpha}_{Y}$-representations, the $L_Y(V_\pi \otimes Y)$'s exhaust the isomorphism classes of simple objects of $\C^G$, since  $$\sum_{Y \in \Irr(\C)/G} \sum_{\pi \in \Irr_{\widetilde{\alpha}_{Y}}} \, \FPdim (L_Y(V_{\pi} \otimes Y))^2 = |G| \, \FPdim \C = \FPdim (\C^G).$$ 	This finishes the proof of the corollary. 
\end{proof}

\medbreak 
Let $Y \in \Irr (\C)/G$ and let $\pi$ be an irreducible
 $\widetilde\alpha_Y$-projective  representation of the group $G_Y$. We shall
use the notation $S_{Y, \pi}$ to indicate the isomorphism class of the simple object of $\C^G$ corresponding to the pair $(Y, \pi)$.
We shall also say that such simple object $S_{Y, \pi}$ \emph{lies
over} $Y$.

\begin{remark}\label{s-object} 
For every set $\mathcal R$ of left coset representatives of $G_Y$
in $G$ and for every collection of isomorphisms $\{c^g: \rho^g(Y)
\to Y\}_{g \in G_Y}$, the class $S_{Y, \pi}$ is represented by the
simple object $S^{\mathcal R, c}_{Y, \pi} : = L_Y^{\mathcal R}(\pi
\otimes Y)$, with the $G_Y$-equivariant structure on $\pi \otimes
Y$  given by \eqref{yip}.

Let us describe more explicitly the dependence of the simple
object $S^{\mathcal R, c}_{Y, \pi}$  on the choice of the
isomorphisms $c^h:\ro^h(Y)\ra Y$. Suppose we are given
another collection of isomorphisms $c'= \{{c'}^g\}$. Then, $Y$
being simple, for any $g \in G_Y$ we can write ${c'}^g = d_{c,
c'}(g) c^g$, for some scalar $d_{c, c'}(g) \in k^*$. It follows
from \eqref{yip} that $\pi \ot Y = d^{-1}_{c, c'}\pi \ot Y$ as
objects of $\cc^{G_Y}$. Hence \bq\label{dyt} S^{\mathcal R,
c}_{Y,\pi} = L^{\mathcal R}_Y(\pi \ot Y) = L^{\mathcal R}_Y(d_{c,
c'}^{-1}\pi \ot Y) = S^{\mathcal R, c'}_{Y, d_{c, c'}^{-1}\pi}.
\eq
\end{remark}

Theorem \ref{correspondence} implies the following:
\begin{lemma}\label{gammaback}
Let $Y \in \Irr(\C)$ and let $\pi$ be an $\widetilde\alpha_Y$-projective representation of $G_Y$. 
Then \beq \pi \simeq \Hom_{\cc}(Y,\;S_{Y, \pi}) \eeq as
 $\widetilde{\alpha}_Y$-projective representations of $G_Y$.
\end{lemma}

\begin{proof} Let $V_\pi$ denote the vector space of the representation $\pi$.
We have $S_{Y, \pi} \simeq L_Y(V_\pi \otimes Y, (\pi(g) \otimes
c^g)_{g \in G})$,  where $c^g:\rho^g(Y)\to Y$ is a collection of
isomorphisms. It follows from Proposition \ref{cociclo} and
Theorem \ref{correspondence} that $\pi \simeq \Hom_\C(Y, V_\pi
\otimes Y) \simeq \Hom_\C(Y, S_{Y, \pi})$.
\end{proof}

As a consequence we now get:

\bp\label{dmultp} Let $Y \in \Irr(\C)$ and let $\pi$ be an irreducible $\alpha_Y$-projective representation of $G_Y$. Then, for all $(X, \mu) \in \C^G$, we have
$$\dim \Hom_{\cc^G}(S_{Y,\pi},\;(X,\mu)) = m_{G_Y}(\pi,\; \Hom_{\cc}(Y, X)).$$
In particular, the simple object $S_{Y, \pi}$ is a constituent of $(X,\mu)$ if and only if $\pi$ is a constituent of $\Hom_{\cc}(Y,X)$.
\ep

Here, $m_{G_Y}(\pi,\; \Hom_{\cc}(Y, X))$ denotes the
multiplicity of $\pi$ in $\Hom_{\cc}(Y, X)$. See Section
\ref{projective-repns}.

\bpf We have a decomposition $(X, \mu) \simeq \oplus_{(Z, \gamma)}\Hom_{\C^G}(S_{Z,\gamma}, (X, \mu)) \otimes S_{Z, \gamma}$, where $Z$ runs over  a set  of representatives of the orbits of the action of $G$ on
$\Irr(\C)$ and $\gamma$ is an irreducible $\alpha_Z$-projective representation of $G_Z$. Since $\Hom_\C(Y, S_{Z,\gamma}) = 0$, for all $Z \neq Y$, then, as projective $G_Y$-representations,
\begin{align*}\Hom_\C(Y, X) & \simeq \oplus_{(Z, \gamma)}\Hom_{\C^G}(S_{Z,\gamma}, (X, \mu)) \otimes \Hom_\C(Y, S_{Z,\gamma})
 \\ & \simeq \oplus_{\gamma}\Hom_{\C^G}(S_{Y,\gamma}, (X, \mu)) \otimes \gamma,
 \end{align*} the last isomorphism by Lemma \ref{gammaback}. This implies the proposition.
\epf

\subsection{On the choice of isomorphisms in a fixed orbit}\label{chosen} Let $Y \in \Irr (\C)$ and let $t \in G$.
Since $\ro^t(Y)$ is a constituent of $F(S_{Y, \pi})$, it follows
from Theorem \ref{correspondence} that  $S_{Y, \pi} \simeq
S_{\ro^t(Y), \delta}$, for some irreducible projective
representation $\delta$ of $G_{\ro^t(Y)}$. In this subsection we
discuss the dependence of $\delta$ upon $\pi$ and the choice of
the sets of isomorphisms $c_Y$, $c_{\ro^t(Y)}$.

\medbreak Let $\pi$ be a projective representation of $G_Y$ with
factor set $\widetilde\al_{Y}$ and let $\;^t\pi$ be the
conjugate projective representation of $G_{\rho^t(Y)} = tG_Yt^{-1}
= : \;^tG_Y$. That is,  $V_{\;^t\pi}=V_{\pi}$ and the action is
defined as $\;^t\pi(h)=\pi(t\inv ht)$, for all $h \in G_Y$. Denote
by $\;^t\widetilde\al_{Y}$ the $2$-cocycle of $G_{\ro^t(Y)}$
given by \bq \;^t\widetilde\al_{Y}(tht^{-1},
th't^{-1})=\widetilde\al_{Y}(h,h'), \quad h,h' \in G_Y. \eq Then
$\;^t\pi$ is a projective representation of $G_{\ro^t(Y)}$ with
factor set $\;^t\widetilde\al_{Y}$.

\medbreak Note that a given collection of isomorphisms
$c^g:\rho^g(Y)\ra Y$, $g \in G_Y$,  determines canonically a
collection of isomorphisms $(^tc)^g:\rho^g(\rho^t(Y))\ra
\rho^t(Y)$, $g \in \;^tG_Y$, in the form \bq\label{ciss} (^tc)^g:=
\rho^t(c^{t^{-1}gt}_Y) \, (\rho_2^{t,t^{-1}gt})^{-1}_Y \,
(\rho_2^{g, t})_Y.\eq Indeed, $t^{-1}gt\in G_Y$ since
$G_{\ro^t(Y)}=tG_Yt^{-1}$.

\begin{remark}\label{alfa-ro} Assume that $Y$ is a simple object representing a class in a fixed orbit
 of the action of $G$. For the objects $\rho^t(Y)$, let the isomorphisms $^tc$ be given as in \eqref{ciss}.
 Then formula \eqref{alfa} gives the $2$-cocycle $\widetilde\al_{\ro^t(Y)}(tht^{-1}, th't^{-1})=\widetilde\al_{Y}(h,h')$ on the inertia subgroup $G_{\ro^t(Y)}$.\end{remark}

\bl \label{xh-equivariant}Let $H$ be a subgroup of $G$ and let
$(M, \nu) \in \C^H$. Then, for all $x \in G$, $(\rho^x(M), ^x\nu)
\in \C^{^xH}$ with equivariant structure $(^x\nu)^{xhx^{-1}}:
\rho^{xhx^{-1}}\rho^x(M) \to \rho^x(M)$ defined, for every $h \in
H$, as the composition \beq
\rho^{xhx^{-1}}(\rho^x(M))\xra{\rho_2^{xhx^{-1}, x}}
\rho^{xh}(M)\xra{(\rho_2^{x,h})^{-1}}\rho^x(\rho^h(M))\xra{\rho^x(\nu^h)}
\rho^x(M). \eeq \el

\begin{proof}
Consider the 2-cocycle $\widetilde\alpha_{\rho^x(M)} \in
Z^2(G_{\rho^x(M)}, k^*)$ associated to the collection of
isomorphisms $^x\nu$. Since $\nu$ is an equivariant structure on
$M$, it has a trivial associated 2-cocycle. It follows from Remark
\ref{alfa-ro} that $\widetilde\alpha_{\rho^x(M)} = 1$ and
therefore $^x\nu$ does define an equivariant structure on
$\rho^x(M)$.
\end{proof}

\begin{corollary} Let $^tc$ be the collection of isomorphisms  given by equation \eqref{ciss} and let $V_{\;^t\pi} \ot \ro^t(Y)$ be the associated object of $\cc^{G_{\ro^t(Y)}}$. Then
$\ro^t(V_{\pi} \ot Y) =  V_{\;^t\pi} \ot \ro^t(Y)$ in $\cc^{G_{\ro^t(Y)}}$.
\end{corollary}

In particular we have an isomorphism of ${}^tG_Y$-equivariant objects $\ro^t(V_{\pi} \ot Y) \simeq V_{\;^t\pi} \ot \ro^t(Y)$, where the ${}^tG_Y$-equivariant structure on $V_{\;^t\pi} \ot \ro^t(Y)$ is induced by any choice of isomorphisms $c_{\rho^t(Y)}$ for $\rho^t(Y)$.

\bpf By Lemma \ref{xh-equivariant}, a $G_{\rho^t(Y)}$-equivariant structure on $\ro^t(V_{\pi} \ot Y)$ is given by
\begin{equation*}
\ro^{tht^{-1}}\rho^t(V_{\pi} \ot Y) \xra{\ro^{tht^{-1},t}_2}\ro^{th}(V_{\pi} \ot Y) \xra{(\rho_2^{t, h})^{-1}} \ro^t \rho^h(V_{\pi} \ot Y) \xra{\pi(h) \otimes \rho^t(c_Y^h)}
V_{\pi} \ot \ro^t(Y),
\end{equation*}
for every $h \in H$. Since $\;^t\pi(tht^{-1})=\pi(h)$, for all $h \in H$, this coincides with the equivariant structure of $V_{\;^t\pi} \ot \ro^t(Y)$ induced by $^tc$.
\epf

%
%
%
%
%
%
%
Suppose now that for every $h \in G_Y$ we have  arbitrary isomorphisms
\begin{equation*}
c_{\ro^t(Y)}^{tht^{-1}}:\ro^{tht^{-1}}(\ro^t(Y))\ra \ro^t(Y)
\end{equation*}
These give rise to isomorphisms
\begin{equation*}
\ro^t(Y)\xra{(c^{tht\inv}_{\ro^t(Y)})\inv} \ro^{tht\inv}(\ro^t(Y))\xra{\ro^{tht\inv , t}_2} \ro^{th}(Y) \xra{(\ro^{t,h}_2)\inv} \ro^t(\ro^h(Y)) \xra{\ro^t(c_Y^h)}\ro^t(Y).
\end{equation*}
Since $\ro^t(Y)$ is a simple object, there exist scalars $d_Y(t,h)\in k^*$ such that
\bq\label{comp}
  \ro^t(c_Y^h) (\ro^{t,h}_2)^{-1} \ro^{tht\inv, t}_2 = d_Y(t,h) c^{tht\inv}_{\ro^t(Y)},
\eq
for all $h \in H$.

\section{Fusion rules for $\cc^G$}\label{fusion-rules}

Let $\pi: G \to \GL(V)$ be a finite dimensional representation of  $G$ on the vector space $V$.
Then the (trivial) object $V \otimes \1 \in \C$ has a $G$-equivariant structure defined by $\pi(g) \otimes \id_\1: \rho^g(V \otimes \1) \to V \otimes \1$.
This induces an embedding of fusion categories $\rep G \to \C^G$ that gives rise to an exact sequence of fusion categories
\begin{equation}\label{sec} \rep G \to \C^G \to \C.
\end{equation}
See \cite[Subsection 5.4]{tensor-exact}.

\begin{remark}\label{invertibles}Let $\Pic(\C)$ be the set of
isomorphism classes of invertible objects of $\C$. The exact
sequence \eqref{sec} induces an exact sequence of groups $$1 \to
\widehat G \to \Pic(\C^G) \to \Pic_0(\C) \to 1,$$ where
$\widehat G \simeq G/[G, G]$ denotes the group of invertible
characters of $G$ and $\Pic_0(\C)$ is the subgroup of $\Pic(\C)$
consisting of isomorphism classes of invertible objects which are
$G$-equivariant. Indeed, $F$ preserves Frobenius-Perron
dimensions, and thus it induces a group homomorphism $F:
\Pic(\C^G) \to \Pic_0(\C)$, which is clearly surjective. The
kernel of $F$ coincides with the invertible objects of $\KER_F =
\rep G$.
\end{remark}

\begin{remark} Note that if $\pi$ is an irreducible representation of $G = G_\1$ on $V$, then the simple object $(V \otimes \1, (\pi(g) \otimes \id_\1)_g)$ of $\C^G$ is isomorphic to
the simple object $S_{\1, \pi}$ corresponding to the pair $(\1, \pi)$ as in Corollary \ref{cor-clifford}.
\end{remark}

\subsection{Orbit formula for the tensor product of two simple objects} Let $Y, Z, U \in \Irr(\C)$ and let $\pi$,
$\gamma$, $\delta$, be projective representations of the
corresponding inertia subgroups with factor sets determined
by \eqref{alfa}.
Let also $S_{Y, \pi}$, $S_{Z,\gm}$ and $S_{U,\delta}$ be the associated simple objects of $\cc^G$.

The multiplicity of $S_{U,\delta}$ in the tensor product
$S_{Y,\pi}\ot S_{Z, \gm}$ is given by the dimension of the vector
space $\Hom_{\cc^G}(S_{U,\delta}, \; S_{Y,\pi}\ot S_{Z,\gm})$. In
view of Proposition \ref{dmultp}, this multiplicity is the same as
the multiplicity of $\delta$ in the space \bq \Hom_{\cc}(U, \;
S_{Y,\pi}\ot S_{Z,\gm}),\eq regarded as an $\al_{U}$-projective
representation of $G_U$.

\medbreak Consider the diagonal action of $G$ on $G/G_Y\times G/G_Z$ coming from the natural actions by left multiplication of $G$ on $G/G_Y$ and $G/G_Z$.  The stabilizer of a pair $(t, s)$ is the subgroup ${}^tG_Y\cap {}^sG_Z \subseteq G$.


\medbreak As objects of $\C$, we have that
\begin{eqnarray*}
S_{Y, \pi} \otimes S_{Z,\gamma}  & \simeq & ( \bigoplus_{u \in G/G_Y}\rho^{u}(V_{\pi}\ot Y))\ot (\bigoplus_{v \in G/G_Z}\rho^{v}(V_{\gm}\ot Z))
\\ & \simeq & \bigoplus_{(u, v) \in G/G_Y\times G/G_Z} \rho^u(V_{\pi}\ot Y) \ot \rho^v(V_{\gm}\ot Z)
\\ & = & \bigoplus_{\cO} S_\cO,
\end{eqnarray*}
where the last summation is over the distinct $G$-orbits $\cO$ in $G/G_Y \times G/G_Z$, and for every $G$-orbit $\cO$,
\bq\label{s-o} S_\cO:=\oplus_{(u, v)\in \cO} \, \rho^u(V_{\pi}\ot Y) \ot \rho^v(V_{\gm}\ot Z).\eq

\medbreak The subgroup $G_U$ acts on $G/G_Y\times G/G_Z$ by restriction and every $G$-orbit in $G/G_Y\times G/G_Z$ is a disjoint union of $G_U$-orbits.
Note that the stabilizer of $(t,s) \in G/G_Y\times G/G_Z$ under the action of $G_U$ is the subgroup $T = G_U \cap tG_Yt^{-1} \cap sG_Zs^{-1}$.

\begin{lemma}\label{mi} For  every $G$-orbit (respectively, $G_U$-orbit)
$\cO \subseteq G/G_Y\times G/G_Z$, $S_\cO$  is an equivariant
subobject (respectively, a $G_U$-equivariant subobject) of the
tensor product $S_{Y, \pi} \otimes S_{Z,\gamma}$.
\end{lemma}

\begin{proof}  We shall prove the statement for $G$-orbits, the proof for $G_U$-orbits being analogous. Let $g \in G$. The equivariant structure $\mu^g: = \mu^g_{ _{S_{Y,\pi}\ot S_{Z,\gamma}}}$ of $S_{Y, \pi} \otimes S_{Z,\gamma}$ is given componentwise by
$$(\mu^g)^{u,v}= (\mu^{g, u}_{S_{Y, \pi}}\ot \mu^{g, v}_{S_{Z,\gamma} }) (V_\pi \otimes V_\gamma \otimes (\rho^g_2)_{\rho^uY, \rho^vZ}),$$ where, for every $(u, v) \in G/G_Y \times G/G_Z$, $\mu^{g, u}_{S_{Y, \pi}}$ and $\mu^{g, v}_{S_{Z,\gamma}}$ are given by formula \eqref{mu}.
It follows that
$$(\mu^g)^{u,v} (\rho^g(V_{\pi}\ot V_{\delta}\ot \rho^u(Y)\ot \rho^v(Z)) \subseteq V_{\pi}\ot V_{\delta}\ot \rho^{u'}(Y)\ot \rho^{v'}(Z),$$ where $(u', v') \in G/G_Y \times G/G_Z$ are uniquely determined by the relations $gu=u'h_Y$ and $gv=v'h_Z$, with $h_Y\in G_Y$ and $h_Z \in G_Z$.
Therefore, $\mu^g (S_\cO) \subseteq S_\cO$, for all $g \in G$.  This implies the lemma. \end{proof}

\medbreak The  map $G \times G \to G$, $(a, b) \mapsto a^{-1}b$, induces a surjective map $p: G/G_Y\times G/G_Z \to G_Y\backslash G/G_Z$, such that $p(tG_Y, sG_Z) = G_Yt^{-1}sG_Z$.

\medbreak Let $\cO_G(t, s)$ denote the $G$-orbit of an element $(t, s) \in G/G_Y\times G/G_Z$.
 Observe that for all $g \in G$, we have $p^{-1}(G_YgG_Z) = \cO_G(e, g)$. Therefore, $p$ induces an identification between the orbit space of $G/G_Y\times G/G_Z$ under the action of $G$ and the space of double cosets $G_Y\backslash G/G_Z$. Combining this with \eqref{s-o}, we obtain:

\begin{corollary}\label{double-cosets} We have a decomposition $S_{Y, \pi} \otimes S_{Z,\gamma} \simeq \bigoplus_{D \in G_Y\backslash G/G_Z} S_D$, where  $S_D : = \bigoplus_{t^{-1}s \in D} \rho^t(V_{\pi}\ot Y) \ot \rho^s(V_{\gm}\ot Z)$.
Moreover, for all $D \in G_Y\backslash G/G_Z$, $S_D$ is a $G$-equivariant subobject of $S_{Y, \pi} \otimes S_{Z,\gamma}$. \qed \end{corollary}

\subsection{Projective representation on multiplicity spaces.}\label{newproj}
It follows from Lemma \ref{mi} that for every $G$-orbit $\cO
\subseteq G/G_Y\times G/G_Z$,   the space $\mathcal H_\cO:=
\Hom_{\cc}(U, S_\cO)$ is an $\al_{U}$-projective representation of
$G_U$.

Let $\cO = \cO_1 \cup \dots \cup \cO_{n_\cO}$ be the decomposition
of $\cO$ into disjoint $G_U$-orbits $\cO_1,  \dots, \cO_{n_\cO}$.
Then, for all $1\leq i \leq n_\cO$, the space $\mathcal H(i) =
\Hom_{\cc}(U, S_{\cO_i})$ is also an $\alpha_U$-projective
representation of $G_U$, where $S_{\cO_i}:=\oplus_{(u, v)\in
\cO_i} \, \rho^u(V_{\pi}\ot Y) \ot \rho^v(V_{\gm}\ot Z)$.

Furthermore, as $G_U$-projective representations, $$\Hom_{\cc}(U,
\; S_{Y, \pi} \otimes S_{Z,\gamma}) \simeq \oplus_{\cO}\mathcal
H_\cO \simeq \oplus_{\cO}\oplus_{i=1}^{n_{\cO}}\mathcal H(i),$$
where summation is understood to run over all orbits $\cO =
\cO_G(t, s)$ such that $\Hom_{\cc}(U, \rho^t(Y)\ot \rho^s(Z)) \neq
0$.

\medbreak For every $(t,s)\in G/G_Y\times G/G_Z$, let  $$\mathcal H_{t,s}:=\Hom_{\cc}(U,\;V_{{}^t\pi} \ot V_{{}^s\gamma}\ot \rho^t(Y)\ot \rho^s(Z)).$$

\begin{lemma}\label{dec-hdei} Let $t, s \in G$. Then $\mathcal H_{t,s}$ is an
$\alpha_U\vert_T$-projective representation of $T = G_U \cap
{}^tG_Y \cap {}^sG_Z$. Moreover, for all $1\leq i\leq n$, we have
\bq\label{h(i)}\mathcal H(i)=\Hom_{\cc}(U,
S_{\cO_i})=\bigoplus_{(t,s)\in \cO_i}\mathcal H_{t,s},\eq as
projective representations of $T$. \end{lemma}

\begin{proof} Since
$V_{{}^t\pi}  \ot \rho^t(Y)$ and $V_{{}^s\gamma}\ot \rho^s(Z)$ are
${}^tG_Y \cap {}^sG_Z$-equivariant objects,  so is their tensor
product. Lemma \ref{proj-rep} implies that $\mathcal H_{t, s}$ is
an $\alpha_U\vert_T$-projective representation of $T$. The
decomposition \eqref{h(i)} follows from the definition of
$S_{\cO_i}$.
\end{proof}


\begin{proposition}\label{new}
Let $U, Y, Z \in \Irr(\C)$ and let $t, s \in G$. Then the vector
space $\tau^{t,s}_U(Y, Z) : = \Hom_{\cc}(U,\rho^t(Y)\ot
\rho^s(Z))$ carries an $\alpha$-projective representation of the
subgroup $T:= G_U\cap G_{\rho^t(Y)} \cap G_{\rho^s(Z)}$, where
$\alpha : =
\al_{U}|_T\al_{\rho^t(Y)}^{-1}|_T\al_{\rho^s(Z)}^{-1}|_T$. The
action of $g \in T$ is given by
\begin{equation}\label{proj-tau}
g.f= (c^g_{\rho^t(Y)} \otimes c^g_{\rho^s(Z)}) \, (\rho^g_2)_{\rho^t(Y),\;\rho^s(Z)} \, \rho^g(f) \, (c^g_U)^{-1},
\end{equation}
for all $f \in \Hom_{\cc}(U,\rho^t(Y)\ot \rho^s(Z))$.
Furthermore, $$\mathcal H_{t, s} \simeq {}^t\pi\vert_T \otimes {}^s\gamma\vert_T \otimes \tau^{t,s}_U(Y, Z),$$ as projective representations of $T$.
\end{proposition}


\begin{remark} Observe that the equivalence class of the projective representation $\tau^{t,s}_U(Y, Z)$ is independent on the choice of the isomorphism classes of $U, Y, Z$ as well as on the choice of isomorphisms $c_{\rho^t(Y)}$, $c_{\rho^s(Z)}$ and $c_U$.
\end{remark}

\bpf Given $X \in \Irr(\C)$, we shall consider in what
follows a fixed (but arbitrary) collection of isomorphisms $c_X =
\{ c^g_X: \rho^g(X) \to X \}_{g \in G_X}$. Let also
$\widetilde\alpha_X \in Z^2(G_X, k^*)$ be the associated
2-cocycles.

We first show that formula \eqref{proj-tau} does define a
projective representation of $T$ with factor set
$\widetilde\al_{U}|_T\widetilde\al_{\rho^t(Y)}^{-1}|_T\widetilde\al_{\rho^s(Z)}^{-1}|_T$.
Let $g, h \in T$, $f \in \Hom_{\cc}(U,\rho^t(Y)\ot \rho^s(Z))$.
Using the definition of the cocycles $\widetilde\alpha$ given by
\eqref{alfa} and relation \eqref{tensor-rho}, we compute:
\begin{align*}g.& (h. f)  = (c^g_{\rho^tY} \otimes c^g_{\rho^sZ}) \, (\rho^g_2)_{\rho^tY,\;\rho^sZ}
\rho^g\left( (c^h_{\rho^tY} \otimes c^h_{\rho^sZ}) \, (\rho^h_2)_{\rho^tY,\;\rho^sZ} \, \rho^h(f) \, {c^h_U}^{-1} \right) {c^g_U}^{-1} \\
& = (c^g_{\rho^tY} \otimes c^g_{\rho^sZ}) \, (\rho^g_2)_{\rho^tY,\;\rho^sZ} \,
\rho^g(c^h_{\rho^tY} \otimes c^h_{\rho^sZ}) \, \rho^g((\rho^h_2)_{\rho^tY,\;\rho^sZ}) \\ & \qquad  \rho^g\rho^h(f) \, \rho^g(c^h_U)^{-1} \, {c^g_U}^{-1} \\
& = (c^g_{\rho^tY} \otimes c^g_{\rho^sZ}) \left(\rho^g(c^h_{\rho^tY}) \otimes \rho^g(c^h_{\rho^sZ})\right)  (\rho^g_2)_{\rho^tY,\;\rho^sZ}
\rho^g((\rho^h_2)_{\rho^tY,\;\rho^sZ}) \, \rho^g\rho^h(f) \\ & \qquad \rho^g(c^h_U)^{-1} {c^g_U}^{-1} \\
& = (c^g_{\rho^tY} \otimes c^g_{\rho^sZ}) \left(\rho^g(c^h_{\rho^tY}) \otimes \rho^g(c^h_{\rho^sZ})\right) \left((\rho_2^{g,h})_{\rho^tY}^{-1} \otimes (\rho_2^{g,h})_{\rho^sZ}^{-1} \right) \, (\rho^{gh}_2)_{\rho^tY,\;\rho^sZ} \\ & \qquad (\rho_2^{g,h})_{\rho^tY \otimes \rho^sZ} \, \rho^g\rho^h(f) \, \rho^g(c^h_U)^{-1} \, {c^g_U}^{-1} \\
& = \widetilde\alpha^{-1}_{\rho^tY}(g, h) \widetilde\alpha^{-1}_{\rho^sZ}(g, h) (c^{gh}_{\rho^tY} \otimes c^{gh}_{\rho^sZ}) \, (\rho^{gh}_2)_{\rho^tY,\;\rho^sZ} \rho^{gh}(f)  (\rho_2^{g,h})_U \, \rho^g(c^h_U)^{-1} \, {c^g_U}^{-1} \\
& = \widetilde\alpha^{-1}_{\rho^tY}(g, h)  \widetilde\alpha^{-1}_{\rho^sZ}(g, h)  \widetilde\alpha_U(g, h) \, (c^{gh}_{\rho^tY} \otimes c^{gh}_{\rho^sZ}) \, (\rho^{gh}_2)_{\rho^tY,\;\rho^sZ} \rho^{gh}(f) \, {c^{gh}_U}^{-1} \\
& = \widetilde\alpha^{-1}_{\rho^tY}(g, h)
\widetilde\alpha^{-1}_{\rho^sZ}(g, h) \widetilde\alpha_U(g, h) \,
(gh. f).
\end{align*}

\medbreak On the other hand, with respect to the given choice of isomorphisms $\{c^g_X\}_{g \in G_X}$, $X\in \Irr (\C)$, the ${}^tG_Y \cap {}^sG_Z$-equivariant structures on $V_{{}^t\pi}  \ot \rho^t(Y)$ and $V_{{}^s\gamma}\ot \rho^s(Z)$ are given, respectively, by ${}^t\pi(g) \otimes c^g_{\rho^t(Y)}: \rho^g(V_{{}^t\pi}  \ot \rho^t(Y)) \to V_{{}^t\pi}  \ot \rho^t(Y)$, and
${}^s\gamma(g) \otimes c^g_{\rho^s(Z)}: \rho^g(V_{{}^s\gamma}  \ot \rho^s(Z)) \to V_{{}^s\gamma}  \ot \rho^s(Z)$, for all $g \in T$.

Thus, the action of $g \in T$ on $f \in \mathcal H_{t, s} =
\Hom_\C(U, V_{{}^t\pi} \otimes V_{{}^s\gamma} \otimes \rho^t(Y)
\otimes \rho^s(Z))$ is determined by
\begin{align*}g.f & =  ({}^t\pi(g) \otimes {}^s\gamma(g) \otimes c^g_{\rho^t(Y)} \otimes c^g_{\rho^s(Z)}) \, (V_{{}^t\pi} \otimes V_{{}^s\gamma} \otimes (\rho^g_2)_{\rho^t(Y), \rho^s(Z)}) \, \rho^g(f) \, (c^g_U)^{-1} \\
& = {}^t\pi(g) \otimes {}^s\gamma(g) \otimes
\left((c^g_{\rho^t(Y)} \otimes c^g_{\rho^s(Z)}) \,
(\rho^g_2)_{\rho^t(Y), \rho^s(Z)}\right) \rho^g(f) \,
(c^g_U)^{-1}.
\end{align*} In view of the $k$-linearity of the functors $\rho^g$, $g \in G$, this implies that the canonical isomorphism $$V_{{}^t\pi} \otimes V_{{}^s\gamma} \otimes \Hom_\C(U, \rho^t(Y) \otimes \rho^s(Z)) \simeq \Hom_\C(U, V_{{}^t\pi} \otimes V_{{}^s\gamma} \otimes  \rho^t(Y) \otimes \rho^s(Z))$$ is indeed an isomorphism of
projective representations of $T$. This finishes the proof of the
proposition. \epf

\begin{proposition}\label{inducida}
Let $(t_i,s_i)\in \cO_i$, $1 \leq i \leq n$, and let $T_i = G_U \cap {}^{t_i}G_Y \cap {}^{s_i}G_Z$ be its stabilizer in $G_U$. Then
$\mathcal H(i)\simeq \Ind_{T_i}^{G_U} \mathcal H_{t_i,s_i}$,
%
as projective $G_U$-representations.
\end{proposition}

\bpf The proposition follows from Lemma \ref{dec-hdei}, in view of
Lemma \ref{perm}. Note that the  group $G_U$ permutes the set
$\cO_i$ transitively. \epf

\subsection{Fusion rules} The following theorem gives the fusion rules for the category $\cc^G$.

\begin{theorem}\label{tp}Let $U, Y, Z \in \Irr (\C)$ and let $\delta, \pi, \gamma$ be irreducible projective representations of the inertia subgroups
$G_U$, $G_Y$, $G_Z$ with factor sets determined by \eqref{alfa}.
Then the multiplicity of $S_{U,\delta}$ in the tensor product $S_{Y, \pi} \otimes S_{Z,\gamma}$ is given by the formula
\begin{equation}\label{multiplicity}\sum_{D \in G_Y\backslash G/G_Z} \sum_{\substack{1\leq i \leq n\\t_i^{-1}s_i \in D \\ \Hom_\C(U, \rho^{t_i}Y \otimes \rho^{s_i}Z) \neq 0}}
m_{T_i}(\delta|_{T_i},\; \;^{t_i}\pi|_{T_i}\ot
\;^{s_i}\gamma|_{T_i} \ot \tau^{s_i,t_i}_U(Y, Z)),\end{equation}
where $(t_1, s_1), \dots, (t_n, s_n)$ are representatives of
the distinct $G_U$-orbits $\cO_1, \dots, \cO_n$ in $G/G_Y \times
G/G_Z$ and, for all $1\leq i \leq n$, $T_i = G_U\cap {}^{t_i}G_Y
\cap {}^{s_i}G_Z$, and $m_{T_i}$ denotes the multiplicity form of
projective $T_i$-representations.
\end{theorem}

\begin{proof} It follows from Proposition \ref{dmultp} that
$$\dim \Hom_{\cc^G}(S_{U,\delta},\;S_{Y, \pi} \otimes S_{Z,\gamma})=
m_{G_U}(\delta,\; \Hom_{\cc}(U, S_{Y, \pi} \otimes
S_{Z,\gamma})).$$ In view of Corollary \ref{double-cosets}, we
have a decomposition
$$\Hom_{\cc}(U, S_{Y, \pi} \otimes S_{Z,\gamma}) \simeq \bigoplus_{D \in G_Y\backslash G/G_Z} \mathcal H_D,$$ as projective representations of $G_U$, where  $$\mathcal H_D : = \bigoplus_{\substack{t^{-1}s \in D, \\ \Hom_\C(U, \rho^t(Y) \ot \rho^s(Z)) \neq 0}} \; \Hom_\C(U, \rho^t(V_{\pi}\ot Y) \ot \rho^s(V_{\gm}\ot Z)).$$

Consider a decomposition $\cO_1 \cup \dots \cup \cO_{n}$ of $G/G_Y \times G/G_Z$ into
disjoint $G_U$-orbits, and  let $\mathcal H(i) \simeq
\Hom_{\cc}(U, S_{\cO_i})$, $1\leq i\leq n$, as in \eqref{h(i)}.

Let also $(t_i, s_i) \in \cO_i$ be a representative of the orbit $\cO_i$ with stabilizer $T_i = G_U\cap {}^{t_i}G_Y \cap {}^{s_i}G_Z$. By Proposition \ref{inducida} we have
$$\mathcal H_D \simeq \bigoplus_{\substack{1\leq i \leq n\\ t_i^{-1}s_i \in D}} \mathcal H(i) \simeq \bigoplus_{\substack{1\leq i \leq n\\ t_i^{-1}s_i \in D}} \Ind_{T_i}^{G_U} \mathcal H_{t_i,s_i},$$
Therefore, using Frobenius Reciprocity and Proposition \ref{new}, we get
\begin{align*}&\dim\Hom_{\cc^G}(S_{U,\delta},  S_{Y, \pi} \otimes S_{Z,\gamma})
= \sum_{D \in G_Y\backslash G/G_Z} m_{G_U}(\delta,\; \mathcal H_D)\\
& = \sum_{D \in G_Y\backslash G/G_Z} \sum_{\substack{1 \leq i \leq n\\ t_i^{-1}s_i \in D \\ \Hom_\C(U, \rho^{t_i}Y \otimes \rho^{s_i}Z) \neq 0}} m_{G_U}(\delta,\; \Ind_{T_i}^{G_U} \mathcal H_{t_i,s_i}) \\
& = \sum_{D \in G_Y\backslash G/G_Z} \sum_{\substack{1 \leq i \leq n\\ t_i^{-1}s_i \in D \\ \Hom_\C(U, \rho^{t_i}Y \otimes \rho^{s_i}Z) \neq 0}} m_{T_i}(\delta\vert_{T_i},\; \mathcal H_{t_i,s_i}) \\
& = \sum_{D \in G_Y\backslash G/G_Z} \sum_{\substack{1 \leq i \leq
n\\ t_i^{-1}s_i \in D \\ \Hom_\C(U, \rho^{t_i}Y \otimes
\rho^{s_i}Z) \neq 0}} m_{T_i}(\delta\vert_{T_i},
{}^{t_i}\pi|_{T_i}\ot {}^{s_i}\gamma|_{T_i} \ot
\tau^{s_i,t_i}_U(Y, Z)). \end{align*}
Thus we get formula \eqref{multiplicity}. This finishes the proof
of the theorem.
\end{proof}

%

\begin{remark}\label{rings} In the paper \cite{ws} a class of rings graded by conjugacy classes of finite groups is studied. More generally, suppose $G$ and $\Gamma$ are finite groups such that $G$ acts on $\Gamma$ by group automorphisms, and denote this action by $g \mapsto {}^xg$, $g \in \Gamma$, $x \in G$. Let $A = \oplus_{g \in \Gamma} A(g)$, where $A(g)$, $g \in \Gamma$, are free modules over a commutative ring $R$, endowed with $R$-linear isomorphisms $c_{g, x}:A(g) \to A({}^xg)$ and $R$-bilinear maps $m_{g, h}: A(g) \times A(h) \to A(gh)$, $g, h \in \Gamma$, $x \in G$, subject to certain compatibility conditions.  Then there is an associative multiplication in $A$ defined componentwise by $m_{g,h}$, $g, h \in G$, and the submodule of $G$-invariants: $A^G = \{ a \in A:\, c_x(a) = a \}$ is a subring of $A$; here $c_x: A \to A$ is the $R$-linear map which is $c_{g,x}$ on the component $A(g)$. 
Furthermore, under weaker assumptions on the maps $c_{g, x}$, $m_{g, h}$,  $A^G$ is an associative ring with multiplication
$$(a.b)_u : = \sum_{\substack{(y, z) \in G_u\backslash (\Gamma \times \Gamma)\\yz = u}} m_{y, z}(a_y, b_z),$$
for all $a, b \in A^G$, $u \in \Gamma$. See \cite[Theorem 2.2]{ws}.

It is shown in \cite[Section 4]{ws} that the Grothendieck ring of the fusion category of representations of a Hopf algebra cocentral abelian extension of $k^\Gamma$ by $kG$ fits into the above construction. 
As shown in \cite[Proposition 3.5]{ext-ty}, this fusion category is an equivariantization of the category $\cc(\Gamma)$ of finite dimensional $\Gamma$-graded vector spaces under the induced action of $G$.

\medbreak Consider now an action of the finite group $G$ on a fusion category $\cc$. For every $Y \in \irr(\cc)$ let $A(Y):=K_0(k_{\al_Y}[G_Y])$ be the Grothendieck group of the twisted group algebra $k_{\al_Y}[G_Y]$ and let $A:=\oplus_{Y \in \irr(\cc)}A(Y)$. For every $g \in G$, let $c_{Y, g}:A(Y) \ra A(\ro^g(Y))$ and $c_g:A\ra A$ be given by
\bq c_{Y, g}(\pi) = {}^g\!\pi, \quad 
c_g:=\oplus_{Y \in G}c_{Y, g}.
\eq 
Theorem \ref{correspondence} implies that there is an isomorphism of abelian groups $\gr(\cc^G) \simeq \{a \in A\:|\; c_Y(a)=a\} = A^G$ such that the simple object $S_{Y, \pi}$ corresponds to the irreducible projective representation $\pi \in A(Y)$, for every $Y \in \Irr(\cc)$.

It follows from Theorem \ref{tp} that, for all $Y, Z \in \Irr(\cc)/G$, the product $S_{Y, \pi}S_{Z, \delta}$ has the expression
\bq\lb{forma}
S_{Y, \pi}S_{Z, \delta} = \sum_{\substack{D\in G_Y \backslash G / G_Z\\ U \in \Irr(\cc)/G}} \sum_{\substack{1 \leq i \leq
n\\ t_i^{-1}s_i \in D \\ \Hom_\C(U, \rho^{t_i}Y \otimes
\rho^{s_i}Z) \neq 0}} S_{U, m^{t_i, s_i}_U(Y, Z)(\pi, \delta)}
\eq
where the map $m^{t,s}_U(Y,Z):A(Y)\times A(Z) \ra  A(U)$ is defined by
\bq
(\pi, \delta)\mapsto \ind^U_T(\;^t\pi|_T\ot \;^s\delta|_T\ot \tau_U^{t,s}(Y,Z))
\eq
Comparing \eqref{forma} with the formula given in \cite[Corollary 2.5]{ws}, we observe that the structure of the Grothendieck ring of the equivariantization $\cc^G$ is similar to that of the rings $A^G$ in \cite{ws}, where the set $\Irr(\cc)$ plays now the r\^ ole of the group $\Gamma$.
\end{remark}

\begin{corollary}\label{constituent} A simple object $S_{U,\delta}$ is a
constituent of a tensor product of simple objects $S_{Y, \pi}
\otimes S_{Z,\gamma}$ if and only if there exist $t\in G/G_Y$ and
$s \in G/G_Z$ such that
\begin{enumerate}\item[(a)] $\Hom_\C(U, \rho^t(Y) \otimes \rho^s(Z)) \neq 0$ and
\item[(b)] $m_{T}(\delta\vert_T,\;^t\pi\vert_T \otimes ^s\gamma\vert_T \otimes \tau^{t, s}_U(Y, Z))\neq 0$, where $T = G_U \cap tG_Yt^{-1} \cap sG_Ys^{-1}$.
\end{enumerate} \qed
\end{corollary}

\subsection{The dual of a simple object}\label{involution}

Let $Y \in \Irr (\C)$. Then the multiplicity of the unit object of
$\cc$ in the tensor product $Y \otimes Y^*$ is one. Hence
$\tau_Y : = \tau^{e, e}_{\1}(Y, Y^*)= \Hom_{\cc}(\1, Y\ot
Y^*)$ is a one dimensional (linear) representation of $G =
G_{\1}$. In particular, it follows from  Proposition \ref{new}
that the cohomology class of the product $\al_{Y} \al_{Y^*}$ is
trivial on $G_Y = G_{Y^*}$.

%
%
%

\medbreak Recall that the dual $\pi^*$ of the $G_Y$-projective
representation  $\pi$ is defined as $V_{\pi}^*$ with $\pi^*(h)(f)
= f \circ \pi(h)^{-1}$. This is an $\al_{Y}^{-1}$-projective
representation of $G_Y$.

\bp\label{dual} The dual object of $S_{Y,\pi}\in \cc^G$ is
determined by 
 $$S_{Y,\pi}^* \simeq S_{Y^*,\;\pi^*\ot\tau_Y^{-1}},$$ where $\tau_Y = \Hom_{\cc}(\1, \;Y\ot Y^*)$. 
\ep

\bpf Observe that $S_{Y,\pi}^* \simeq S_{Z, \gamma}$, for some $Z
\in \Irr(\C)/G$ and some $\alpha_Z$-projective representation of
$G_Z$. On the other hand, $S_{Y,\pi}^* \simeq S_{Z, \gamma}$ if
and only if the unit object is a constituent of $S_{Y,\pi} \otimes
S_{Z, \gamma}$. Since the unit object of $\C^G$ is isomorphic to
$S_{\1, \epsilon}$, where $\epsilon$ denotes the trivial
representation of $G_\1 = G$, it follows from Corollary
\ref{constituent} that $S_{Y,\pi}^* \simeq
S_{Y^*,\;\pi^*\ot\tau_Y^{-1}}$, where $\tau_Y = \Hom_{\cc}(\1,
\;Y\ot Y^*)$, as claimed.  \epf

Combining Proposition \ref{dual} with Frobenius Reciprocity we
obtain:

\begin{corollary}\label{onecst} Let $Y, Z \in \Irr(\C)/G$ and let $\pi$, $\gm$ be projective representations of
$G_Y$ and $G_Z$ with factor sets determined by \eqref{alfa}.
Let also $\delta$ be an irreducible representation
of $G$. Then $S_{\1,\delta}$ is a constituent of $S_{Y, \pi}\ot
S^*_{Z, \gm}$ if and only if $Z = Y^*$ and $\delta$ is a
constituent of $(\pi \ot\gm^*\ot\tau_Z^{-1})\uw^G_{G_Y}$. \qed
\end{corollary}

\subsection{$\C^G$ as a $\rep G$-bimodule
category}\label{bimodcat}

Let us regard the category $\rep G$ as a fusion subcategory of $\C^G$
via the natural embedding $\pi \mapsto (\pi \otimes \1,
\pi(g)\otimes \id_\1)$. So that the tensor product of $\C^G$ makes
$\C^G$ into a $\rep G$-bimodule category.

\medbreak The results in Section \ref{simples-equiv} imply that
there is an equivalence of $k$-linear categories $\C^G \simeq
\oplus_{Y \in \Irr(\C)/G} \, \rep_{\alpha_Y} G_Y$, where
$\rep_{\alpha_Y} G_Y$ is the category of finite dimensional
$\alpha_Y$-projective representations of $G_Y$. Under this
equivalence, a simple object $\pi$ of $\rep_{\alpha_Y} G_Y$, that
is, an irreducible $\alpha_Y$-projective representation of $G_Y$,
corresponds to the simple object $S_{Y, \pi}$ of $\C^G$. In other
words, $\rep_{\alpha_Y} G_Y$ is identified with the full
subcategory of $\C$ whose simple objects are lying over $Y$.
An explicit equivalence is determined, for every $Y \in
\Irr(\C)/G$, by the functors $L_Y: \C^{G_Y} \to \C^G$ and $F_Y:
\C^G \to \C^{G_Y}$.

For each $Y \in \Irr(\C)$, the category $\rep_{\alpha_Y} G_Y$ is
in a canonical way an indecomposable $\rep G$-bimodule category
via tensor product of projective representations; see
\cite[Theorem 3.2]{ostrik}. As a consequence of Theorem \ref{tp}
we obtain:

\begin{theorem}\label{bimod-equiv} There is an equivalence of $\rep G$-bimodule categories \begin{equation}\label{k-lin} \C^G \simeq \bigoplus_{Y \in
\Irr(\C)/G} \, \rep_{\alpha_Y} G_Y.
\end{equation}
Moreover, each $\rep_{\alpha_Y} G_Y$ is an indecomposable $\rep G$-bimodule category.
\end{theorem}

\begin{proof} Let $\pi$ be an irreducible representation of $G$, so
that $\pi$ corresponds to the simple object $S_{\1, \pi} \in \rep
G$, and let $S_{Z, \gamma} \in \rep_{\alpha_Z}G_Z$ be another
simple object, where $Z \in \Irr(\C)/G$. It follows from Corollary
\ref{constituent} that if the simple object $S_{U, \delta}$, $U\in
\Irr(\C)/G$, is a constituent of $S_{\1, \pi} \otimes S_{Z,
\gamma}$, then $U \simeq \rho^s(Z)$, for some $s \in G/G_Z$. Hence
$U = Z$ and thus the group $T = G_Z \cap G_\1 \cap G_Z$ coincides
with $G_Z$,  $\tau_U(\1, Z) \simeq \Hom_\C(Z, Z)$ is a one
dimensional (linear) representation of $G_Z$. Therefore
$\pi\vert_{G_Z} \otimes \gamma \otimes \tau_U(\1, Z) \simeq
\pi\vert_{G_Z} \otimes \gamma$ as projective representations of
$G_Z$.

By Theorem \ref{tp}, the multiplicity of $S_{U, \delta}$ in the
tensor product $S_{\1, \pi} \otimes S_{Z, \gamma}$ equals
$m_{G_Z}(\delta, \pi \otimes \gamma)$. Therefore we obtain
$$S_{\1, \pi} \otimes S_{Z, \gamma} \simeq
\bigoplus_{\delta}\, m_{G_Z}(\delta, \pi\vert_{G_Z} \otimes
\gamma) \; S_{Z, \delta},$$ where $\delta$ runs over the
equivalence classes of $\alpha_Z$-projective representations of
$G_Z$. Clearly this object corresponds to $\pi\vert_{G_Z} \otimes
\gamma \in \rep_{\alpha_Z}G_Z$.

Similar arguments apply for the tensor product $S_{Z, \gamma}
\otimes S_{\1, \pi}$. This completes the proof of the theorem.
\end{proof}

For any $U\in \Irr(\C)$, let us extend the notation $S_{U, \delta}
= L_U(\delta \otimes U)$ to indicate the object of $\C^G$
corresponding to an arbitrary $\alpha_U$-projective representation
$\delta$ of $G_U$.

\begin{remark} Let $Y, Z \in \Irr(\C)$ and let $S_{Y,\pi}$, $S_{Z, \gamma}$ be
simple objects of $\C^G$ lying over $Y$ and $Z$, respectively. So
that $\pi$ is an irreducible $\widetilde\alpha_Y$-projective representation
of $G_Y$ and $\gamma$ is an irreducible $\widetilde\alpha_Z$-projective
representation of $G_Z$.

According to Theorem \ref{bimod-equiv}, the tensor product
$S_{Y,\pi}\otimes S_{Z, \gamma}$ has a decomposition
\begin{equation}\label{dec-bimod}S_{Y,\pi}\otimes S_{Z, \gamma} \cong \bigoplus_{U \in
\Irr(\C)/G}S_{U, \delta},
\end{equation}
where, for all $U \in \Irr(\C)/G$, $S_{U, \delta} \in \C^G$ is the
sum of simple constituents of $S_{Y,\pi}\otimes S_{Z, \gamma}$
lying over $U$. It follows from Proposition \ref{dmultp} that
$\delta \cong \Hom_\C(U, S_{Y,\pi}\otimes S_{Z, \gamma})$. \end{remark}

\begin{remark}
The
action of $G$ on $\C$ induces an action of $G$ on $\gr(\C)$ by
algebra automorphisms. Let $\gr(\C)^G \subseteq \gr(\C)$ be the
subring of $G$-invariants in $\gr(\C)$. For every $Y \in \Irr
(\C)$, let us  consider the element \bq\label{def-csy} \cs(Y):= \sum_{t \in
G/G_Y}\ro^t(Y) \in \gr(\C). \eq Clearly, we have $\cs(Y) \in
\gr(\C)^G$ and $\cs(Y) = \cs(\rho^g(Y))$, for all $Y \in
\Irr(\C)$. Observe that $F^!(S_{Y,\pi})=(\dim \pi )S(Y)$, where $F^!:
\gr(\C^G) \to \gr(\C)$ is the ring map induced by the forgetful
functor $F:\C^G \to \C$.
Moreover, the set $\{ \cs(Y): \, Y \in \Irr(\C)/G \}$
is a basis for $\gr(\C)^G$ and,  for all $Y, Z \in
\Irr(\C)/G$, we have \bq \cs(Y)\cs(Z)=\sum_{U\in \Irr(\C)/G}m_{Y, Z}^U \;
\cs(U), \eq for some nonnegative integers $m^U_{Y,Z}$.

Let $Y, Z, U \in \Irr(\C)/G$. Consider any fixed simple objects $S_{Y, \pi}$ and $S_{Z, \gamma}$ of $\C^G$ lying over $Y$ and $Z$, respectively.
Applying the map $F^!$ in formula \eqref{dec-bimod}, we obtain that $m^U_{Y, Z} = \dim \delta / (\dim \pi) (\dim \gamma)$, where
$\delta = \Hom_\C(U, S_{Y,\pi}\otimes S_{Z, \gamma}) \simeq V_\pi \otimes V_\gamma \otimes (\oplus_{(t, s) \in G/G_Y \times G/G_Z}\Hom_\C(U, \rho^t(Y) \otimes \rho^s(Z)))$.
Therefore, for all $Y, Z, U \in \Irr(\C)/G$, the integers $m^U_{Y, Z}$ are given
by the formula $$m^U_{Y, Z} = \sum_{(t, s) \in G/G_Y \times G/G_Z} \dim \Hom_\C(U, \rho^t(Y) \otimes \rho^s(Z)).$$
\end{remark}

\section{Application to equivariantizations of pointed fusion categories}\label{pointed-fusion}

We shall consider in this section a \emph{pointed} fusion category $\C$, that is,  all simple objects of $\C$ are invertible. Then there is an equivalence of fusion categories $\C \simeq \C(\Gamma, \omega)$, where $\Gamma = \Pic(\C)$ is the group of isomorphism classes of invertible objects in $\C$,   $\omega: \Gamma \times \Gamma \times \Gamma \to k^*$ is an invertible normalized 3-cocycle and $\C(\Gamma, \omega) = \vect^\Gamma_\omega$ is the category of finite dimensional $\Gamma$-graded vector spaces with associativity constraint induced by $\omega$.

\subsection{Group actions on $\C(\Gamma, \omega)$ and equivariantizations}

Let $\C = \C(\Gamma, \omega)$ and let $G$ be a finite group. An action $\rho: \underline G \to \underline \Aut_{\otimes} \C$ of $G$ on $\C$ is determined by an action by group automorphisms of $G$ on $\Gamma$,  that we shall indicate by $x \mapsto {}^gx$, $x \in \Gamma$, $g \in G$, and two maps $\tau: G \times \Gamma \times \Gamma \to k^*$, and $\sigma: G \times G \times \Gamma \to k^*$, satisfying
\begin{align*}\frac{\omega(x, y, z)}{\omega({}^gx, {}^gy, {}^gz)} & = \frac{\tau(g; xy, z) \, \tau(g; x, y)}{\tau(g; y, z) \, \tau(g; x, yz)} \\
1 & = \frac{\sigma(h, l; x) \, \sigma(g, hl; x)}{\sigma(gh, l; x) \, \sigma(g, h; {}^lx)} \\
\frac{\tau(gh; x, y)}{\tau(g; {}^hx, {}^hy) \, \tau(h; x, y)} & = \frac{\sigma(g, h; x) \, \sigma(g, h; y)}{\sigma(g, h; xy)},
\end{align*} for all $x, y, z \in \Gamma$, $g, h, l \in G$.

We shall also assume that $\tau$ and $\sigma$ satisfy the additional normalization conditions $\tau(g; x, y) = \sigma(g, h; x) = 1$, whenever some of the arguments $g$, $h$, $x$ or $y$ is an identity.

\medbreak The action $\rho: \underline G \to \underline \Aut_{\otimes} \C$ determined by this data is defined by letting $\rho^g(x) = {}^gx$, for all $g\in G$, $x\in \Gamma$, and $\rho^g = \id$ on arrows, together with the following constraints:
\begin{equation}\label{formulas-pointed}(\rho_2^{g, h})_x = \sigma(g, h; x)^{-1} \id_{^{gh}x}, \quad (\rho_2^g)_{x, y} = \tau(g; x, y)^{-1}\id_{xy}, \quad \rho^g_0 = \id_e,
\end{equation} for all $g, h \in G$, $x, y \in \Gamma$. See \cite[Section 7]{tambara}.

\subsection{Fusion rules for $\C(\Gamma, \omega)^G$}
Let us denote $\sigma_x(g, h): = \sigma(g, h; x)$ and $\tau_{x, y}(g): = \tau(g; x, y)$, $x, y \in \Gamma$, $g, h \in G$.

\medbreak For all $x \in \Gamma$ and $g \in G_x$ we let the
isomorphism $c_x: {}^gx = x \to x$  to be the identity of $x$.
Therefore, the cocycle $\widetilde\alpha_x: G_x \times G_x
\to k^*$ defined by \eqref{alfa} is given by \bq\label{al-x}
\widetilde\alpha_x(g, h) = \sigma_x(g, h)^{-1},\eq for all $g, h
\in G_x$.

\medbreak It follows from Corollary \ref{cor-clifford} that the set of isomorphism classes of simple objects of $\C^G$ is parameterized by isomorphism classes of pairs $(y, \pi)$, where $y$ runs over the orbits of the action of $G$ on $\Gamma$ and $\pi$ is an irreducible projective representation of the inertia subgroup $G_y \subseteq G$ with factor set $\sigma_y$.

\medbreak Let $\cO$ be a  $G$-orbit in $G/G_y\times G/G_z$ corresponding to a double coset $D \in G_y\backslash G/G_z$. Then $\cO = \cO_G(e, g)$, for any $g \in D$, and $\cO$ contains at most one $G_U$-orbit, $\cO_{G_U}(t, s)$, $t^{-1}s \in D$, such that $\Hom_\C(x, {}^ty \otimes {}^sz) \neq 0$. Indeed, the condition $\Hom_\C(x, {}^ty \otimes {}^sz) \neq 0$ amounts in this case to $x =  {}^ty{}^sz$. Thus, for all $e\neq g \in G/G_U$, $x \neq {}^gx = {}^{gt}y{}^{gs}z$.

\medbreak In addition, if $x = {}^{t}y{}^{s}z$, then ${}^{t}G_y \cap {}^{s}G_z \subseteq G_x$. Therefore, $G_x\cap {}^{t}G_y \cap {}^{s}G_z = {}^{t}G_y \cap {}^{s}G_z$.

\medbreak In the projective representation of ${}^{t}G_y \cap {}^{s}G_z$ in $\Hom_\C(x, {}^{t}y{}^{s}z) \simeq k$, defined in Lemma \ref{new}, the action of an element $g \in {}^{t}y{}^{s}z$ is nothing but scalar multiplication by $\tau_{{}^ty, {}^sz}(g)^{-1}$.

\medbreak As a consequence of  Theorem \ref{tp}, the following theorem gives the fusion rules for the category $\C(\Gamma, \omega)^G$.

\begin{theorem}\label{tp-pointed}Let $x, y, z \in \Gamma$ and let $\delta, \pi, \gamma$ be
irreducible projective representations of the inertia subgroups
$G_x$, $G_y$, $G_z$ with factor sets $\sigma_x$, $\sigma_y$,
$\sigma_z$, respectively. Then the multiplicity of $S_{x,\delta}$
in the tensor product $S_{y, \pi} \otimes S_{z,\gamma}$ is given
by the formula
\begin{equation*}\sum_{D \in G_y\backslash G/G_z} \sum_{\substack{t^{-1}s \in D \\ x = {}^{t}y{}^{s}z}}
m_{\, {}^{t}G_y \cap {}^{s}G_z}(\delta|_{{}^{t}G_y \cap {}^{s}G_z},\; \;^{t}\pi|_{{}^{t}G_y \cap {}^{s}G_z}\ot \;^{s}\gamma|_{{}^{t}G_y \cap {}^{s}G_z} \; \tau_{{}^ty, {}^sz}^{-1}).\end{equation*}
\end{theorem}

\begin{example} Consider a cocentral abelian exact sequence of Hopf algebras $k \to k^\Gamma \to H \to kG \to k$, where $\Gamma$ and $F$ are finite groups.
As special case of  \cite[Proposition 3.5]{ext-ty}, there is an
action of $G$ on the category $\C = \C(\Gamma, 1)$ of finite
dimensional representations of $k^\Gamma$ such that $\Rep H \simeq
\C^G$ as fusion categories (see \cite[Remark 2.1]{ext-ty}). In
this situation, the formula for the fusion rules of $\Rep H$ given
by Theorem \ref{tp-pointed} specializes to the formula obtained by
C. Goff in \cite[Theorem 4.5]{goff}. See Section \ref{projective-repns}. \end{example}

\subsection{Braided group-theoretical fusion categories}
Recall that a fusion category is called \emph{group-theoretical}
if it is Morita equivalent to a pointed fusion category. In view of \cite[Theorem 7.2]{NNW}, a \emph{braided}
fusion category is group-theoretical if and only  if it is an
equivariantization of a pointed fusion category.

More precisely, it  was shown in \cite[Theorem 5.3]{naidu}  that
every braided group-theoretical fusion category is equivalent to
an equivariantization $\C(\xi)^G$ of a crossed pointed fusion
category $\C(\xi)$ associated to a quasi-abelian 3-cocycle $\xi$
on a finite crossed module $(G, X, \partial)$, under a canonical
action of $G$ on $\C(\xi)$.

\medbreak Recall that a finite crossed module $(G, X, \partial)$ consists of a finite group $G$ acting by automorphisms on a finite group $X$, and a group homomorphism $\partial: X \to G$ such that $${}^{\partial(x)}y = xyx^{-1}, \quad \partial({}^gx) = g\partial(x)g^{-1}, \qquad g \in G,\; x, y \in X,$$ where $x \mapsto ^gx$, $x \in X$, $g \in G$, denotes the action of $g$ on $X$.

A quasi-abelian 3-cocycle $\xi$ on $(G, X, \partial)$ is a quadruple $\xi = (\omega, \gamma, \mu, c)$, where $\omega: X \times X \times X \to k^*$ is a $3$-cocycle, $\gamma: G \times G \times X \to k^*$, $\mu: G \times X \times X \to k^*$ and $c: X \times X \to k^*$ are maps satisfying the compatibility conditions in \cite[Definition 3.4]{naidu}.

\medbreak As a fusion category $\C(\xi) = \C(X, \omega)$, and the action of $G$ on $\C(\xi)$ is determined by the action of $G$ on $X$ and formulas \eqref{formulas-pointed}, with respect to $\sigma_x(g, h) : = \gamma(g, h; x)$, $\tau_{x, y}(g) : = \mu(g; x, y)^{-1}$, $x, y \in X$, $g, h \in G$. See \cite[Subsection 4.1]{naidu}.

\medbreak Theorem \ref{tp-pointed} gives thus the fusion rules in
the category $\C(\xi)^G$ in terms of group-theoretical data
determined by the crossed module $(G, X, \partial)$ and the
quasi-abelian 3-cocycle $\xi$, entailing the determination of the
fusion rules in any braided group-theoretical fusion category.

\begin{example} Let $\omega: G \times G \times G \to k^*$ be a 3-cocycle on $G$.
Consider the crossed module $(G, G, \id)$ with respect to the
adjoint action of $G$ on itself. The quadruple $\xi =
(\omega^{-1}, \gamma^{-1}, \mu^{-1}, 1)$ is a quasi-abelian
3-cocycle on $(G, G, \id)$, where $\gamma$ and $\mu$ are defined
in the form
\begin{align*}\gamma(g, h; x) & = \frac{\omega(g, h, x) \omega(ghxh^{-1}g^{-1}, g, h)}{\omega(g, hxh^{-1}, y)},
\\ \mu(g; x, y) & = \frac{\omega(gxg^{-1}, g,
y)}{\omega(gxg^{-1}, gyg^{-1}, g)\omega(g, x, y)},
\end{align*} for all $g, h, x, y \in G$.

The equivariantization $\C(\xi)^G$ is equivalent to the category
$\Rep D^{\omega}G$ of finite dimensional representations of the
twisted quantum double $D^{\omega}G$ introduced in \cite{dpr}. See
\cite[Lemma 6.3]{naidu}.

\medbreak Simple objects of $\C(\xi)^G$ are parameterized by
$S_{x, \pi}$, where $x$ runs over a set of representatives of
conjugacy classes of $G$ and $\pi$ is an irreducible projective
representation of the centralizer $Z(x)$ of $x$ in $G$ with
factor set $\gamma_x$. Theorem \ref{tp-pointed} gives the following
formula for the multiplicity of $S_{x, \delta}$ in the tensor
product $S_{y, \pi} \otimes S_{z, \gamma}$:
\begin{equation*}\sum_{D \in Z(y)\backslash G/Z(z)} \sum_{\substack{t^{-1}s \in D \\ x = tyt^{-1} szs^{-1}}}
m_{\, tZ(y)t^{-1} \cap sZ(z)s^{-1}}(\delta,\; \;^{t}\pi \ot
\;^{s}\gamma \; \mu(-; tyt^{-1}, szs^{-1})^{-1}).\end{equation*}

We point out that the fusion rules for the category $\Rep
D^{\omega}G$ were also determined in \cite[Section 5]{goff}.
\end{example}

\section{Appendix}\label{projective-repns}
In this Appendix we give a brief account of the results on projective representations used in the paper.
See for instance \cite{karpilovsky}.

\medbreak Let $G$ be a finite group and let
$\widetilde\alpha: G \times G \to k^*$ be a (normalized)
2-cocycle on $G$, that is,
$$\widetilde\alpha(g, h)\widetilde\alpha(gh, t) = \widetilde\alpha(g, ht) \widetilde\alpha(h, t), \quad \widetilde\alpha(g, e) = 1 = \widetilde\alpha(e, g),\qquad \forall g, h, t \in G.$$

\medbreak A \emph{projective representation} $\pi$ of $G$ with
\emph{factor set} $\widetilde\alpha$
 on a vector space $V$ is a map $\pi: G\ra
\GL(V)$, such that $$\pi(e) = \id_V, \quad
\pi(gh)=\widetilde\al(g,\;h)\pi(g)\pi(h),\qquad \forall g, h\in
G.$$
In other words, $\pi$ is a representation of the twisted group
algebra $k_{\widetilde\alpha} G$ on the vector space $V$. We shall
also use the notation  $V_{\pi} = V$ to indicate such a projective
representation.  Thus two projective $\widetilde{\alpha}$-representations  $\pi$ and $\pi'$ are isomorphic if and only if they are isomorphic $k_{\widetilde\alpha} G$-representations. In
this case we shall use the notation $\pi \simeq \pi'$. 

\medbreak Two projective representations $\pi$ and $\pi'$ of $G$
are called \emph{(projectively) equivalent}  if there is a linear
isomorphism $\phi: V_\pi \to V_{\pi'}$ and a map $f: G \to k^*$
such that $\phi \pi(g) = f(g)\pi'(g) \phi$, for all $g \in G$. 

If  $\pi'$ and $\pi$ are projectively equivalent, then the associated cocycles
$\widetilde\alpha$ and $\widetilde\alpha'$ are related by
$$\widetilde\alpha(g, h) = \widetilde\alpha'(g, h)f(g)f(h)f(gh)^{-1}, \quad g, h \in
G,$$  that is, $\widetilde\alpha$ and $\widetilde\alpha'$ are
cohomologous cocycles, and thus they belong to the same cohomology
class $\alpha \in H^2(G, k^*)$. We shall also call $\pi$ an
$\alpha$-projective representation. Note that the map $f: G \to
k^*$ induces an algebra isomorphism $\tilde f:
k_{\widetilde\alpha} G \to k_{\widetilde\alpha'} G$ in the form
$\tilde f(g) = f(g)g$, for all $g \in G$. Thus $\pi$ and $\pi'$
are equivalent projective representations if and only if $V_{\pi}
\simeq {\tilde f}^*(V_{\pi'})$ as $k_{\widetilde\alpha}G$-modules.

\medbreak Let $\pi$ and $\pi'$ be projective representations of
$G$ with factor sets $\widetilde\alpha$ and $\widetilde\alpha'$,
respectively. The \emph{tensor product} $\pi \otimes \pi'$ is the
projective $\widetilde\alpha \widetilde\alpha'$-representation on
the vector space $V_{\pi} \otimes V_{\pi'}$ defined by $(\pi
\otimes \pi')(g) (u \otimes v) = \pi(g)u \otimes \pi'(g)v$. In
particular, if $\pi$ is a representation of $G$, then $\pi \otimes
\pi'$ is again a projective representation with factor set
$\widetilde\alpha$.

If $\pi_1$ and $\pi_1'$ are projective representations
projectively equivalent to $\pi$ and $\pi'$, respectively, then
the tensor products $\pi_1\otimes \pi'_1$ and $\pi\otimes \pi'$
are projectively equivalent. Further, suppose that $\pi'$ is a
one-dimensional representation, that is, a linear character of
$G$. Then  $\pi$ and  $\pi \otimes \pi'$ are projectively
equivalent via the canonical isomorphism $\phi: V_{\pi} \to
V_{\pi} \otimes k$, $v \mapsto v \otimes 1$, and the map $f: G \to
k^*$ given by $f(g) = \pi'(g)^{-1}$, for all $g \in G$.

\medbreak A nonzero projective representation $\pi: G \to \GL(V)$
of $G$ is called \emph{irreducible} if $0$ and $V$ are the only
subspaces of $V$ which are invariant under $\pi(g)$, for all $g
\in G$. Hence, $\pi$ is irreducible if and only if it is not
projectively equivalent to a projective representation
$\rho$ of the form $$\rho(g) = \left(\begin{matrix}\pi_1(g) & * \\
0 & \pi_2(g)
\end{matrix}\right), \quad g \in G,$$ where $\pi_1$ and $\pi_2$
are nonzero projective representations or, equivalently, if $V$ is
a simple $k_{\widetilde\alpha} G$-module, where $\widetilde
\alpha$ is the factor set of $\pi$ \cite[Theorem
3.2.5]{karpilovsky}.

\medbreak Let $\pi: G \to \GL(V)$ be a projective representation
of $G$ with factor set $\widetilde\alpha$. Since the group algebra
$k_{\widetilde\alpha} G$ is semisimple, then $V = V_\pi$ is
\emph{completely reducible}, that is, $V_\pi \simeq V_{\pi_1}
\oplus \dots \oplus V_{\pi_n}$, where $V_{\pi_i}$ is a simple
$k_{\widetilde\alpha} G$-module, for all $i = 1, \dots, n$. If
$\pi'$ is an irreducible projective representation with factor set
$\widetilde\alpha'$, then $\pi'$ is called a \emph{constituent} of
$\pi$ if $\pi'$ is projectively equivalent to $\pi_i$ for some
$1\leq i \leq n$. In this case, the \emph{multiplicity} (or
\emph{intertwining number}) of $\pi'$ in $\pi$ is defined as
$$m_G(\pi', \pi) : = \dim \Hom_{k_{\widetilde\alpha} G}(V_{\pi_i}, V_{\pi}).$$
Observe that if $\pi'$ is a constituent of $\pi$, then the
cocycles $\widetilde\alpha'$ and  $\widetilde\alpha$ belong to the
same class  in $H^2(G, k^*)$. Letting $\widetilde\alpha' df =
\widetilde\alpha$, with $f: G \to k^*$, we have that $m_G(\pi',
\pi) : = \dim \Hom_{k_{\widetilde\alpha} G}({\tilde
f}^*(V_{\pi'}), V_{\pi})$, where $\tilde f: k_{\widetilde\alpha} G
\to k_{\widetilde\alpha'} G$ is the isomorphism associated  to
$f$.

\medbreak The \emph{character} of a projective representation
$\pi: G \to \GL(V)$ is defined as the map $\chi = \chi_V: G \to
k^*$ given by $\chi(g) = \Tr(\pi(g))$, for all $g \in G$. Let
$\widetilde\alpha$ be the factor set of $\pi$. If $\pi'$ is an
irreducible projective representation of $G$ with factor set
$\widetilde\alpha$ and character $\chi'$, then the multiplicity of
$\pi'$ in $\pi$ can be computed by the formula
\begin{align*}m_G(\pi', \pi)  = \langle \chi', \chi \rangle : & = \frac{1}{|G|} \sum_{g \in G}  \frac{1}{\widetilde\alpha(g^{-1}, g)} \chi'(g) \chi(g^{-1})  \\
& = \frac{1}{|G|} \sum_{g \in G^0}
\frac{1}{\widetilde\alpha(g^{-1}, g)} \chi'(g) \chi(g^{-1}),
\end{align*} $G^0 \subseteq G$ is the subset of $\widetilde\alpha$-regular elements of $G$. See \cite[Chapter 5]{karpilovsky}.

\medbreak Let $\widetilde\alpha: G \times G \to k^*$ be a
2-cocycle and let $H \subseteq G$ be a subgroup. Consider a
projective representation $W$ of $H$ with factor set
$\widetilde\alpha\vert_H$. The \emph{induced projective
representation} of $G$ is defined as $\Ind_H^GW =
k_{\widetilde\alpha}G \otimes_{k_{\widetilde\alpha} H}W$. This is
a projective representation of $G$ with factor set
$\widetilde\alpha$. By Frobenius reciprocity, we have  natural
isomorphisms
\begin{equation*} \Hom_{k_{\widetilde\alpha} G}(\Ind_H^GW, V) \simeq \Hom_{k_{\widetilde\alpha} H}(W, V\vert_H),\end{equation*}
for every projective representation $V$ of $G$ with factor set
$\widetilde\alpha$, where $V\vert_H$ denotes the restricted
projective representation of $H$.

The following lemma gives a characterization of those projective representations which are induced from a subgroup.

\begin{lemma}\label{perm} Let $\widetilde\alpha: G \times G \to k^*$ be a cocycle and let $V$ be a $k_{\widetilde\alpha} G$-module.
Suppose   $V = \oplus_{x \in X}V_x$ is a grading of $V$ by a set
$X$ and assume that there is transitive action of $G$ on $X$, $G
\times X \to X$, $(g, x) \mapsto {}^gx$, such that $g. V_x =
V_{^gx}$, for all $g \in G$, $x \in X$.

Let also $y \in X$ and $G_y \subseteq G$ the inertia subgroup of
$y$. Then $V_y$ is a $k_{\widetilde\alpha} G_y$-module and $V
\simeq \Ind_{G_y}^GV_y$ as $k_{\widetilde\alpha} G$-modules.
\end{lemma}

\bpf See \cite[Theorem 5.2.1]{karpilovsky}.
%
%
\epf


\begin{thebibliography}{A}

\bibitem{agaitsgory} {\sc S. Arkhipov} and {\sc D. Gaitsgory},
\emph{Another realization of the category of modules over the
small quantum group}, Adv. Math. \textbf{173}, 114--143 (2003).

\bibitem{tensor-exact} {\sc A. Brugui\`{e}res} and {\sc S. Natale}, \emph{Exact sequences of tensor categories},
Int. Math. Res. Not. \textbf{2011}  (24), 5644--5705 (2011).

\bibitem{indp-exact} {\sc A. Brugui\`{e}res} and {\sc S. Natale}, \emph{Central exact sequences of tensor categories, equivariantization and applications},
J. Math. Soc. Japan, to appear. Preprint    \texttt{arXiv:1112.3135}.


\bibitem{dpr}  {\sc R. Dijkgraaf} and {\sc V. Pasquier} and {\sc Ph. Roche, } \emph{Quasi-Quantum Groups Related to Orbifold Models}, Proceedings of the ``International Colloquium on Modern Quantum Field Theory,'' Tata Institute of Fundamental Research, 375--383 (1990).

\bibitem{DGNOI}  {\sc V. Drinfeld}, {\sc S. Gelaki} , {\sc D. Nikshych} and {\sc V. Ostrik},
\emph{On braided fusion categories I}, Sel. Math. New Ser.
\textbf{16}, 1--119 (2010).

\bibitem{ENO2}  {\sc P. Etingof}, {\sc D. Nikshych} and {\sc V. Ostrik},
\emph{Weakly group-theoretical and solvable fusion categories}, Adv. Math. \textbf{226}, 176--205 (2011).

\bibitem{fw}  {\sc E. Frenkel}, and {\sc E. Witten},
\emph{Geometric endoscopy and mirror symmetry}, Commun. Number
Theory Phys. \textbf{2}, 113--283 (2008).


\bibitem{goff} {\sc C. Goff}, \emph{Fusion rules for abelian extensions of Hopf algebras}, Algebra Number Theory \textbf{6},  327--339 (2012).


\bibitem{karpilovsky} {\sc G. Karpilovsky}, \emph{Proyective Representation of Finite Groups},
Pure and Applied Mathematics \textbf{94}, Marcel Dekker, New
York-Basel (1985).

\bibitem{NNW} {\sc D. Naidu}, \textsc{D. Nikshych} and {\sc S. Witherspoon}, \emph{Fusion subcategories of representation categories of
twisted quantum doubles of finite groups}, Int. Math. Res. Not. \textbf{2009}, 4183--4219 (2009).

\bibitem{naidu} {\sc D. Naidu}, \emph{Crossed pointed categories and their equivariantizations}, Pacific J.  Math. \textbf{247}, 477-496 (2010).


\bibitem{ext-ty} {\sc S. Natale}, \emph{Hopf algebra extensions of group algebras and Tambara-Yamagami categories}, Algebr. Represent. Theory \textbf{13},  673--691 (2010).

\bibitem{nik} \textsc{D. Nikshych},
\emph{Non group-theoretical semisimple Hopf algebras from group
actions on fusion categories}, Sel. Math. New Ser. \textbf{14}, 145--161
(2008).

\bibitem{MW} {\sc S. Montgomery} and {\sc S. Witherspoon}, \emph{Irreducible Representations of Crossed Products}, J. Pure and Appl. Algebra \textbf{111}, 381--385 (1988).

\bibitem{ostrik} {\sc V. Ostrik}, \emph{Module categories, weak Hopf algebras and modular invariants},
Transform. Groups \textbf{8}, 177--206 (2003).

\bibitem{tambara}  \textsc{D. Tambara},
\emph{Invariants and semi-direct products for finite group actions
on tensor categories}, J. Math. Soc. Japan \textbf{53},
429--456 (2001).

\bibitem{ws} {\sc S. Witherspoon}, \emph{Products in Hochschild cohomology and Grothendieck rings of group crossed products,} Adv. Math. \textbf{185}, 136--158 (2004).


\end{thebibliography}
\end{document}